\documentclass[bj,preprint]{imsart}

\RequirePackage[OT1]{fontenc}
\RequirePackage{amsthm,amsmath,natbib}
\RequirePackage[colorlinks,citecolor=blue,urlcolor=blue]{hyperref}

\usepackage{color}  
\definecolor{AliceBlue}{rgb}{0.94,0.97,1.00}
\definecolor{AntiqueWhite1}{rgb}{1.00,0.94,0.86}
\definecolor{AntiqueWhite2}{rgb}{0.93,0.87,0.80}
\definecolor{AntiqueWhite3}{rgb}{0.80,0.75,0.69}
\definecolor{AntiqueWhite4}{rgb}{0.55,0.51,0.47}
\definecolor{AntiqueWhite}{rgb}{0.98,0.92,0.84}
\definecolor{BlanchedAlmond}{rgb}{1.00,0.92,0.80}
\definecolor{BlueViolet}{rgb}{0.54,0.17,0.89}
\definecolor{CadetBlue1}{rgb}{0.60,0.96,1.00}
\definecolor{CadetBlue2}{rgb}{0.56,0.90,0.93}
\definecolor{CadetBlue3}{rgb}{0.48,0.77,0.80}
\definecolor{CadetBlue4}{rgb}{0.33,0.53,0.55}
\definecolor{CadetBlue}{rgb}{0.37,0.62,0.63}
\definecolor{CornflowerBlue}{rgb}{0.39,0.58,0.93}
\definecolor{DarkBlue}{rgb}{0.00,0.00,0.55}
\definecolor{DarkCyan}{rgb}{0.00,0.55,0.55}
\definecolor{DarkGoldenrod1}{rgb}{1.00,0.73,0.06}
\definecolor{DarkGoldenrod2}{rgb}{0.93,0.68,0.05}
\definecolor{DarkGoldenrod3}{rgb}{0.80,0.58,0.05}
\definecolor{DarkGoldenrod4}{rgb}{0.55,0.40,0.03}
\definecolor{DarkGoldenrod}{rgb}{0.72,0.53,0.04}
\definecolor{DarkGray}{rgb}{0.66,0.66,0.66}
\definecolor{DarkGreen}{rgb}{0.00,0.39,0.00}
\definecolor{DarkGrey}{rgb}{0.66,0.66,0.66}
\definecolor{DarkKhaki}{rgb}{0.74,0.72,0.42}
\definecolor{DarkMagenta}{rgb}{0.55,0.00,0.55}
\definecolor{DarkOliveGreen1}{rgb}{0.79,1.00,0.44}
\definecolor{DarkOliveGreen2}{rgb}{0.74,0.93,0.41}
\definecolor{DarkOliveGreen3}{rgb}{0.64,0.80,0.35}
\definecolor{DarkOliveGreen4}{rgb}{0.43,0.55,0.24}
\definecolor{DarkOliveGreen}{rgb}{0.33,0.42,0.18}
\definecolor{DarkOrange1}{rgb}{1.00,0.50,0.00}
\definecolor{DarkOrange2}{rgb}{0.93,0.46,0.00}
\definecolor{DarkOrange3}{rgb}{0.80,0.40,0.00}
\definecolor{DarkOrange4}{rgb}{0.55,0.27,0.00}
\definecolor{DarkOrange}{rgb}{1.00,0.55,0.00}
\definecolor{DarkOrchid1}{rgb}{0.75,0.24,1.00}
\definecolor{DarkOrchid2}{rgb}{0.70,0.23,0.93}
\definecolor{DarkOrchid3}{rgb}{0.60,0.20,0.80}
\definecolor{DarkOrchid4}{rgb}{0.41,0.13,0.55}
\definecolor{DarkOrchid}{rgb}{0.60,0.20,0.80}
\definecolor{DarkRed}{rgb}{0.55,0.00,0.00}
\definecolor{DarkSalmon}{rgb}{0.91,0.59,0.48}
\definecolor{DarkSeaGreen1}{rgb}{0.76,1.00,0.76}
\definecolor{DarkSeaGreen2}{rgb}{0.71,0.93,0.71}
\definecolor{DarkSeaGreen3}{rgb}{0.61,0.80,0.61}
\definecolor{DarkSeaGreen4}{rgb}{0.41,0.55,0.41}
\definecolor{DarkSeaGreen}{rgb}{0.56,0.74,0.56}
\definecolor{DarkSlateBlue}{rgb}{0.28,0.24,0.55}
\definecolor{DarkSlateGray1}{rgb}{0.59,1.00,1.00}
\definecolor{DarkSlateGray2}{rgb}{0.55,0.93,0.93}
\definecolor{DarkSlateGray3}{rgb}{0.47,0.80,0.80}
\definecolor{DarkSlateGray4}{rgb}{0.32,0.55,0.55}
\definecolor{DarkSlateGray}{rgb}{0.18,0.31,0.31}
\definecolor{DarkSlateGrey}{rgb}{0.18,0.31,0.31}
\definecolor{DarkTurquoise}{rgb}{0.00,0.81,0.82}
\definecolor{DarkViolet}{rgb}{0.58,0.00,0.83}
\definecolor{DeepPink1}{rgb}{1.00,0.08,0.58}
\definecolor{DeepPink2}{rgb}{0.93,0.07,0.54}
\definecolor{DeepPink3}{rgb}{0.80,0.06,0.46}
\definecolor{DeepPink4}{rgb}{0.55,0.04,0.31}
\definecolor{DeepPink}{rgb}{1.00,0.08,0.58}
\definecolor{DeepSkyBlue1}{rgb}{0.00,0.75,1.00}
\definecolor{DeepSkyBlue2}{rgb}{0.00,0.70,0.93}
\definecolor{DeepSkyBlue3}{rgb}{0.00,0.60,0.80}
\definecolor{DeepSkyBlue4}{rgb}{0.00,0.41,0.55}
\definecolor{DeepSkyBlue}{rgb}{0.00,0.75,1.00}
\definecolor{DimGray}{rgb}{0.41,0.41,0.41}
\definecolor{DimGrey}{rgb}{0.41,0.41,0.41}
\definecolor{DodgerBlue1}{rgb}{0.12,0.56,1.00}
\definecolor{DodgerBlue2}{rgb}{0.11,0.53,0.93}
\definecolor{DodgerBlue3}{rgb}{0.09,0.45,0.80}
\definecolor{DodgerBlue4}{rgb}{0.06,0.31,0.55}
\definecolor{DodgerBlue}{rgb}{0.12,0.56,1.00}
\definecolor{FloralWhite}{rgb}{1.00,0.98,0.94}
\definecolor{ForestGreen}{rgb}{0.13,0.55,0.13}
\definecolor{GhostWhite}{rgb}{0.97,0.97,1.00}
\definecolor{GreenYellow}{rgb}{0.68,1.00,0.18}
\definecolor{HotPink1}{rgb}{1.00,0.43,0.71}
\definecolor{HotPink2}{rgb}{0.93,0.42,0.65}
\definecolor{HotPink3}{rgb}{0.80,0.38,0.56}
\definecolor{HotPink4}{rgb}{0.55,0.23,0.38}
\definecolor{HotPink}{rgb}{1.00,0.41,0.71}
\definecolor{IndianRed1}{rgb}{1.00,0.42,0.42}
\definecolor{IndianRed2}{rgb}{0.93,0.39,0.39}
\definecolor{IndianRed3}{rgb}{0.80,0.33,0.33}
\definecolor{IndianRed4}{rgb}{0.55,0.23,0.23}
\definecolor{IndianRed}{rgb}{0.80,0.36,0.36}
\definecolor{LavenderBlush1}{rgb}{1.00,0.94,0.96}
\definecolor{LavenderBlush2}{rgb}{0.93,0.88,0.90}
\definecolor{LavenderBlush3}{rgb}{0.80,0.76,0.77}
\definecolor{LavenderBlush4}{rgb}{0.55,0.51,0.53}
\definecolor{LavenderBlush}{rgb}{1.00,0.94,0.96}
\definecolor{LawnGreen}{rgb}{0.49,0.99,0.00}
\definecolor{LemonChiffon1}{rgb}{1.00,0.98,0.80}
\definecolor{LemonChiffon2}{rgb}{0.93,0.91,0.75}
\definecolor{LemonChiffon3}{rgb}{0.80,0.79,0.65}
\definecolor{LemonChiffon4}{rgb}{0.55,0.54,0.44}
\definecolor{LemonChiffon}{rgb}{1.00,0.98,0.80}
\definecolor{LightBlue1}{rgb}{0.75,0.94,1.00}
\definecolor{LightBlue2}{rgb}{0.70,0.87,0.93}
\definecolor{LightBlue3}{rgb}{0.60,0.75,0.80}
\definecolor{LightBlue4}{rgb}{0.41,0.51,0.55}
\definecolor{LightBlue}{rgb}{0.68,0.85,0.90}
\definecolor{LightCoral}{rgb}{0.94,0.50,0.50}
\definecolor{LightCyan1}{rgb}{0.88,1.00,1.00}
\definecolor{LightCyan2}{rgb}{0.82,0.93,0.93}
\definecolor{LightCyan3}{rgb}{0.71,0.80,0.80}
\definecolor{LightCyan4}{rgb}{0.48,0.55,0.55}
\definecolor{LightCyan}{rgb}{0.88,1.00,1.00}
\definecolor{LightGoldenrod1}{rgb}{1.00,0.93,0.55}
\definecolor{LightGoldenrod2}{rgb}{0.93,0.86,0.51}
\definecolor{LightGoldenrod3}{rgb}{0.80,0.75,0.44}
\definecolor{LightGoldenrod4}{rgb}{0.55,0.51,0.30}
\definecolor{LightGoldenrodYellow}{rgb}{0.98,0.98,0.82}
\definecolor{LightGoldenrod}{rgb}{0.93,0.87,0.51}
\definecolor{LightGray}{rgb}{0.83,0.83,0.83}
\definecolor{LightGreen}{rgb}{0.56,0.93,0.56}
\definecolor{LightGrey}{rgb}{0.83,0.83,0.83}
\definecolor{LightPink1}{rgb}{1.00,0.68,0.73}
\definecolor{LightPink2}{rgb}{0.93,0.64,0.68}
\definecolor{LightPink3}{rgb}{0.80,0.55,0.58}
\definecolor{LightPink4}{rgb}{0.55,0.37,0.40}
\definecolor{LightPink}{rgb}{1.00,0.71,0.76}
\definecolor{LightSalmon1}{rgb}{1.00,0.63,0.48}
\definecolor{LightSalmon2}{rgb}{0.93,0.58,0.45}
\definecolor{LightSalmon3}{rgb}{0.80,0.51,0.38}
\definecolor{LightSalmon4}{rgb}{0.55,0.34,0.26}
\definecolor{LightSalmon}{rgb}{1.00,0.63,0.48}
\definecolor{LightSeaGreen}{rgb}{0.13,0.70,0.67}
\definecolor{LightSkyBlue1}{rgb}{0.69,0.89,1.00}
\definecolor{LightSkyBlue2}{rgb}{0.64,0.83,0.93}
\definecolor{LightSkyBlue3}{rgb}{0.55,0.71,0.80}
\definecolor{LightSkyBlue4}{rgb}{0.38,0.48,0.55}
\definecolor{LightSkyBlue}{rgb}{0.53,0.81,0.98}
\definecolor{LightSlateBlue}{rgb}{0.52,0.44,1.00}
\definecolor{LightSlateGray}{rgb}{0.47,0.53,0.60}
\definecolor{LightSlateGrey}{rgb}{0.47,0.53,0.60}
\definecolor{LightSteelBlue1}{rgb}{0.79,0.88,1.00}
\definecolor{LightSteelBlue2}{rgb}{0.74,0.82,0.93}
\definecolor{LightSteelBlue3}{rgb}{0.64,0.71,0.80}
\definecolor{LightSteelBlue4}{rgb}{0.43,0.48,0.55}
\definecolor{LightSteelBlue}{rgb}{0.69,0.77,0.87}
\definecolor{LightYellow1}{rgb}{1.00,1.00,0.88}
\definecolor{LightYellow2}{rgb}{0.93,0.93,0.82}
\definecolor{LightYellow3}{rgb}{0.80,0.80,0.71}
\definecolor{LightYellow4}{rgb}{0.55,0.55,0.48}
\definecolor{LightYellow}{rgb}{1.00,1.00,0.88}
\definecolor{LimeGreen}{rgb}{0.20,0.80,0.20}
\definecolor{MediumAquamarine}{rgb}{0.40,0.80,0.67}
\definecolor{MediumBlue}{rgb}{0.00,0.00,0.80}
\definecolor{MediumOrchid1}{rgb}{0.88,0.40,1.00}
\definecolor{MediumOrchid2}{rgb}{0.82,0.37,0.93}
\definecolor{MediumOrchid3}{rgb}{0.71,0.32,0.80}
\definecolor{MediumOrchid4}{rgb}{0.48,0.22,0.55}
\definecolor{MediumOrchid}{rgb}{0.73,0.33,0.83}
\definecolor{MediumPurple1}{rgb}{0.67,0.51,1.00}
\definecolor{MediumPurple2}{rgb}{0.62,0.47,0.93}
\definecolor{MediumPurple3}{rgb}{0.54,0.41,0.80}
\definecolor{MediumPurple4}{rgb}{0.36,0.28,0.55}
\definecolor{MediumPurple}{rgb}{0.58,0.44,0.86}
\definecolor{MediumSeaGreen}{rgb}{0.24,0.70,0.44}
\definecolor{MediumSlateBlue}{rgb}{0.48,0.41,0.93}
\definecolor{MediumSpringGreen}{rgb}{0.00,0.98,0.60}
\definecolor{MediumTurquoise}{rgb}{0.28,0.82,0.80}
\definecolor{MediumVioletRed}{rgb}{0.78,0.08,0.52}
\definecolor{MidnightBlue}{rgb}{0.10,0.10,0.44}
\definecolor{MintCream}{rgb}{0.96,1.00,0.98}
\definecolor{MistyRose1}{rgb}{1.00,0.89,0.88}
\definecolor{MistyRose2}{rgb}{0.93,0.84,0.82}
\definecolor{MistyRose3}{rgb}{0.80,0.72,0.71}
\definecolor{MistyRose4}{rgb}{0.55,0.49,0.48}
\definecolor{MistyRose}{rgb}{1.00,0.89,0.88}
\definecolor{NavajoWhite1}{rgb}{1.00,0.87,0.68}
\definecolor{NavajoWhite2}{rgb}{0.93,0.81,0.63}
\definecolor{NavajoWhite3}{rgb}{0.80,0.70,0.55}
\definecolor{NavajoWhite4}{rgb}{0.55,0.47,0.37}
\definecolor{NavajoWhite}{rgb}{1.00,0.87,0.68}
\definecolor{NavyBlue}{rgb}{0.00,0.00,0.50}
\definecolor{OldLace}{rgb}{0.99,0.96,0.90}
\definecolor{OliveDrab1}{rgb}{0.75,1.00,0.24}
\definecolor{OliveDrab2}{rgb}{0.70,0.93,0.23}
\definecolor{OliveDrab3}{rgb}{0.60,0.80,0.20}
\definecolor{OliveDrab4}{rgb}{0.41,0.55,0.13}
\definecolor{OliveDrab}{rgb}{0.42,0.56,0.14}
\definecolor{OrangeRed1}{rgb}{1.00,0.27,0.00}
\definecolor{OrangeRed2}{rgb}{0.93,0.25,0.00}
\definecolor{OrangeRed3}{rgb}{0.80,0.22,0.00}
\definecolor{OrangeRed4}{rgb}{0.55,0.15,0.00}
\definecolor{OrangeRed}{rgb}{1.00,0.27,0.00}
\definecolor{PaleGoldenrod}{rgb}{0.93,0.91,0.67}
\definecolor{PaleGreen1}{rgb}{0.60,1.00,0.60}
\definecolor{PaleGreen2}{rgb}{0.56,0.93,0.56}
\definecolor{PaleGreen3}{rgb}{0.49,0.80,0.49}
\definecolor{PaleGreen4}{rgb}{0.33,0.55,0.33}
\definecolor{PaleGreen}{rgb}{0.60,0.98,0.60}
\definecolor{PaleTurquoise1}{rgb}{0.73,1.00,1.00}
\definecolor{PaleTurquoise2}{rgb}{0.68,0.93,0.93}
\definecolor{PaleTurquoise3}{rgb}{0.59,0.80,0.80}
\definecolor{PaleTurquoise4}{rgb}{0.40,0.55,0.55}
\definecolor{PaleTurquoise}{rgb}{0.69,0.93,0.93}
\definecolor{PaleVioletRed1}{rgb}{1.00,0.51,0.67}
\definecolor{PaleVioletRed2}{rgb}{0.93,0.47,0.62}
\definecolor{PaleVioletRed3}{rgb}{0.80,0.41,0.54}
\definecolor{PaleVioletRed4}{rgb}{0.55,0.28,0.36}
\definecolor{PaleVioletRed}{rgb}{0.86,0.44,0.58}
\definecolor{PapayaWhip}{rgb}{1.00,0.94,0.84}
\definecolor{PeachPuff1}{rgb}{1.00,0.85,0.73}
\definecolor{PeachPuff2}{rgb}{0.93,0.80,0.68}
\definecolor{PeachPuff3}{rgb}{0.80,0.69,0.58}
\definecolor{PeachPuff4}{rgb}{0.55,0.47,0.40}
\definecolor{PeachPuff}{rgb}{1.00,0.85,0.73}
\definecolor{PowderBlue}{rgb}{0.69,0.88,0.90}
\definecolor{RosyBrown1}{rgb}{1.00,0.76,0.76}
\definecolor{RosyBrown2}{rgb}{0.93,0.71,0.71}
\definecolor{RosyBrown3}{rgb}{0.80,0.61,0.61}
\definecolor{RosyBrown4}{rgb}{0.55,0.41,0.41}
\definecolor{RosyBrown}{rgb}{0.74,0.56,0.56}
\definecolor{RoyalBlue1}{rgb}{0.28,0.46,1.00}
\definecolor{RoyalBlue2}{rgb}{0.26,0.43,0.93}
\definecolor{RoyalBlue3}{rgb}{0.23,0.37,0.80}
\definecolor{RoyalBlue4}{rgb}{0.15,0.25,0.55}
\definecolor{RoyalBlue}{rgb}{0.25,0.41,0.88}
\definecolor{SaddleBrown}{rgb}{0.55,0.27,0.07}
\definecolor{SandyBrown}{rgb}{0.96,0.64,0.38}
\definecolor{SeaGreen1}{rgb}{0.33,1.00,0.62}
\definecolor{SeaGreen2}{rgb}{0.31,0.93,0.58}
\definecolor{SeaGreen3}{rgb}{0.26,0.80,0.50}
\definecolor{SeaGreen4}{rgb}{0.18,0.55,0.34}
\definecolor{SeaGreen}{rgb}{0.18,0.55,0.34}
\definecolor{SkyBlue1}{rgb}{0.53,0.81,1.00}
\definecolor{SkyBlue2}{rgb}{0.49,0.75,0.93}
\definecolor{SkyBlue3}{rgb}{0.42,0.65,0.80}
\definecolor{SkyBlue4}{rgb}{0.29,0.44,0.55}
\definecolor{SkyBlue}{rgb}{0.53,0.81,0.92}
\definecolor{SlateBlue1}{rgb}{0.51,0.44,1.00}
\definecolor{SlateBlue2}{rgb}{0.48,0.40,0.93}
\definecolor{SlateBlue3}{rgb}{0.41,0.35,0.80}
\definecolor{SlateBlue4}{rgb}{0.28,0.24,0.55}
\definecolor{SlateBlue}{rgb}{0.42,0.35,0.80}
\definecolor{SlateGray1}{rgb}{0.78,0.89,1.00}
\definecolor{SlateGray2}{rgb}{0.73,0.83,0.93}
\definecolor{SlateGray3}{rgb}{0.62,0.71,0.80}
\definecolor{SlateGray4}{rgb}{0.42,0.48,0.55}
\definecolor{SlateGray}{rgb}{0.44,0.50,0.56}
\definecolor{SlateGrey}{rgb}{0.44,0.50,0.56}
\definecolor{SpringGreen1}{rgb}{0.00,1.00,0.50}
\definecolor{SpringGreen2}{rgb}{0.00,0.93,0.46}
\definecolor{SpringGreen3}{rgb}{0.00,0.80,0.40}
\definecolor{SpringGreen4}{rgb}{0.00,0.55,0.27}
\definecolor{SpringGreen}{rgb}{0.00,1.00,0.50}
\definecolor{SteelBlue1}{rgb}{0.39,0.72,1.00}
\definecolor{SteelBlue2}{rgb}{0.36,0.67,0.93}
\definecolor{SteelBlue3}{rgb}{0.31,0.58,0.80}
\definecolor{SteelBlue4}{rgb}{0.21,0.39,0.55}
\definecolor{SteelBlue}{rgb}{0.27,0.51,0.71}
\definecolor{VioletRed1}{rgb}{1.00,0.24,0.59}
\definecolor{VioletRed2}{rgb}{0.93,0.23,0.55}
\definecolor{VioletRed3}{rgb}{0.80,0.20,0.47}
\definecolor{VioletRed4}{rgb}{0.55,0.13,0.32}
\definecolor{VioletRed}{rgb}{0.82,0.13,0.56}
\definecolor{WhiteSmoke}{rgb}{0.96,0.96,0.96}
\definecolor{YellowGreen}{rgb}{0.60,0.80,0.20}
\definecolor{aliceblue}{rgb}{0.94,0.97,1.00}
\definecolor{antiquewhite}{rgb}{0.98,0.92,0.84}
\definecolor{aquamarine1}{rgb}{0.50,1.00,0.83}
\definecolor{aquamarine2}{rgb}{0.46,0.93,0.78}
\definecolor{aquamarine3}{rgb}{0.40,0.80,0.67}
\definecolor{aquamarine4}{rgb}{0.27,0.55,0.45}
\definecolor{aquamarine}{rgb}{0.50,1.00,0.83}
\definecolor{azure1}{rgb}{0.94,1.00,1.00}
\definecolor{azure2}{rgb}{0.88,0.93,0.93}
\definecolor{azure3}{rgb}{0.76,0.80,0.80}
\definecolor{azure4}{rgb}{0.51,0.55,0.55}
\definecolor{azure}{rgb}{0.94,1.00,1.00}
\definecolor{beige}{rgb}{0.96,0.96,0.86}
\definecolor{bisque1}{rgb}{1.00,0.89,0.77}
\definecolor{bisque2}{rgb}{0.93,0.84,0.72}
\definecolor{bisque3}{rgb}{0.80,0.72,0.62}
\definecolor{bisque4}{rgb}{0.55,0.49,0.42}
\definecolor{bisque}{rgb}{1.00,0.89,0.77}
\definecolor{black}{rgb}{0.00,0.00,0.00}
\definecolor{blanchedalmond}{rgb}{1.00,0.92,0.80}
\definecolor{blue1}{rgb}{0.00,0.00,1.00}
\definecolor{blue2}{rgb}{0.00,0.00,0.93}
\definecolor{blue3}{rgb}{0.00,0.00,0.80}
\definecolor{blue4}{rgb}{0.00,0.00,0.55}
\definecolor{blueviolet}{rgb}{0.54,0.17,0.89}
\definecolor{blue}{rgb}{0.00,0.00,1.00}
\definecolor{brown1}{rgb}{1.00,0.25,0.25}
\definecolor{brown2}{rgb}{0.93,0.23,0.23}
\definecolor{brown3}{rgb}{0.80,0.20,0.20}
\definecolor{brown4}{rgb}{0.55,0.14,0.14}
\definecolor{brown}{rgb}{0.65,0.16,0.16}
\definecolor{burlywood1}{rgb}{1.00,0.83,0.61}
\definecolor{burlywood2}{rgb}{0.93,0.77,0.57}
\definecolor{burlywood3}{rgb}{0.80,0.67,0.49}
\definecolor{burlywood4}{rgb}{0.55,0.45,0.33}
\definecolor{burlywood}{rgb}{0.87,0.72,0.53}
\definecolor{cadetblue}{rgb}{0.37,0.62,0.63}
\definecolor{chartreuse1}{rgb}{0.50,1.00,0.00}
\definecolor{chartreuse2}{rgb}{0.46,0.93,0.00}
\definecolor{chartreuse3}{rgb}{0.40,0.80,0.00}
\definecolor{chartreuse4}{rgb}{0.27,0.55,0.00}
\definecolor{chartreuse}{rgb}{0.50,1.00,0.00}
\definecolor{chocolate1}{rgb}{1.00,0.50,0.14}
\definecolor{chocolate2}{rgb}{0.93,0.46,0.13}
\definecolor{chocolate3}{rgb}{0.80,0.40,0.11}
\definecolor{chocolate4}{rgb}{0.55,0.27,0.07}
\definecolor{chocolate}{rgb}{0.82,0.41,0.12}
\definecolor{coral1}{rgb}{1.00,0.45,0.34}
\definecolor{coral2}{rgb}{0.93,0.42,0.31}
\definecolor{coral3}{rgb}{0.80,0.36,0.27}
\definecolor{coral4}{rgb}{0.55,0.24,0.18}
\definecolor{coral}{rgb}{1.00,0.50,0.31}
\definecolor{cornflowerblue}{rgb}{0.39,0.58,0.93}
\definecolor{cornsilk1}{rgb}{1.00,0.97,0.86}
\definecolor{cornsilk2}{rgb}{0.93,0.91,0.80}
\definecolor{cornsilk3}{rgb}{0.80,0.78,0.69}
\definecolor{cornsilk4}{rgb}{0.55,0.53,0.47}
\definecolor{cornsilk}{rgb}{1.00,0.97,0.86}
\definecolor{cyan1}{rgb}{0.00,1.00,1.00}
\definecolor{cyan2}{rgb}{0.00,0.93,0.93}
\definecolor{cyan3}{rgb}{0.00,0.80,0.80}
\definecolor{cyan4}{rgb}{0.00,0.55,0.55}
\definecolor{cyan}{rgb}{0.00,1.00,1.00}
\definecolor{darkblue}{rgb}{0.00,0.00,0.55}
\definecolor{darkcyan}{rgb}{0.00,0.55,0.55}
\definecolor{darkgoldenrod}{rgb}{0.72,0.53,0.04}
\definecolor{darkgray}{rgb}{0.66,0.66,0.66}
\definecolor{darkgreen}{rgb}{0.00,0.39,0.00}
\definecolor{darkgrey}{rgb}{0.66,0.66,0.66}
\definecolor{darkkhaki}{rgb}{0.74,0.72,0.42}
\definecolor{darkmagenta}{rgb}{0.55,0.00,0.55}
\definecolor{darkolive}{rgb}{0.33,0.42,0.18}
\definecolor{darkorange}{rgb}{1.00,0.55,0.00}
\definecolor{darkorchid}{rgb}{0.60,0.20,0.80}
\definecolor{darkred}{rgb}{0.55,0.00,0.00}
\definecolor{darksalmon}{rgb}{0.91,0.59,0.48}
\definecolor{darksea}{rgb}{0.56,0.74,0.56}
\definecolor{darkslate}{rgb}{0.18,0.31,0.31}
\definecolor{darkslate}{rgb}{0.18,0.31,0.31}
\definecolor{darkslate}{rgb}{0.28,0.24,0.55}
\definecolor{darkturquoise}{rgb}{0.00,0.81,0.82}
\definecolor{darkviolet}{rgb}{0.58,0.00,0.83}
\definecolor{deeppink}{rgb}{1.00,0.08,0.58}
\definecolor{deepsky}{rgb}{0.00,0.75,1.00}
\definecolor{dimgray}{rgb}{0.41,0.41,0.41}
\definecolor{dimgrey}{rgb}{0.41,0.41,0.41}
\definecolor{dodgerblue}{rgb}{0.12,0.56,1.00}
\definecolor{firebrick1}{rgb}{1.00,0.19,0.19}
\definecolor{firebrick2}{rgb}{0.93,0.17,0.17}
\definecolor{firebrick3}{rgb}{0.80,0.15,0.15}
\definecolor{firebrick4}{rgb}{0.55,0.10,0.10}
\definecolor{firebrick}{rgb}{0.70,0.13,0.13}
\definecolor{floralwhite}{rgb}{1.00,0.98,0.94}
\definecolor{forestgreen}{rgb}{0.13,0.55,0.13}
\definecolor{gainsboro}{rgb}{0.86,0.86,0.86}
\definecolor{ghostwhite}{rgb}{0.97,0.97,1.00}
\definecolor{gold1}{rgb}{1.00,0.84,0.00}
\definecolor{gold2}{rgb}{0.93,0.79,0.00}
\definecolor{gold3}{rgb}{0.80,0.68,0.00}
\definecolor{gold4}{rgb}{0.55,0.46,0.00}
\definecolor{goldenrod1}{rgb}{1.00,0.76,0.15}
\definecolor{goldenrod2}{rgb}{0.93,0.71,0.13}
\definecolor{goldenrod3}{rgb}{0.80,0.61,0.11}
\definecolor{goldenrod4}{rgb}{0.55,0.41,0.08}
\definecolor{goldenrod}{rgb}{0.85,0.65,0.13}
\definecolor{gold}{rgb}{1.00,0.84,0.00}
\definecolor{gray0}{rgb}{0.00,0.00,0.00}
\definecolor{gray100}{rgb}{1.00,1.00,1.00}
\definecolor{gray10}{rgb}{0.10,0.10,0.10}
\definecolor{gray11}{rgb}{0.11,0.11,0.11}
\definecolor{gray12}{rgb}{0.12,0.12,0.12}
\definecolor{gray13}{rgb}{0.13,0.13,0.13}
\definecolor{gray14}{rgb}{0.14,0.14,0.14}
\definecolor{gray15}{rgb}{0.15,0.15,0.15}
\definecolor{gray16}{rgb}{0.16,0.16,0.16}
\definecolor{gray17}{rgb}{0.17,0.17,0.17}
\definecolor{gray18}{rgb}{0.18,0.18,0.18}
\definecolor{gray19}{rgb}{0.19,0.19,0.19}
\definecolor{gray1}{rgb}{0.01,0.01,0.01}
\definecolor{gray20}{rgb}{0.20,0.20,0.20}
\definecolor{gray21}{rgb}{0.21,0.21,0.21}
\definecolor{gray22}{rgb}{0.22,0.22,0.22}
\definecolor{gray23}{rgb}{0.23,0.23,0.23}
\definecolor{gray24}{rgb}{0.24,0.24,0.24}
\definecolor{gray25}{rgb}{0.25,0.25,0.25}
\definecolor{gray26}{rgb}{0.26,0.26,0.26}
\definecolor{gray27}{rgb}{0.27,0.27,0.27}
\definecolor{gray28}{rgb}{0.28,0.28,0.28}
\definecolor{gray29}{rgb}{0.29,0.29,0.29}
\definecolor{gray2}{rgb}{0.02,0.02,0.02}
\definecolor{gray30}{rgb}{0.30,0.30,0.30}
\definecolor{gray31}{rgb}{0.31,0.31,0.31}
\definecolor{gray32}{rgb}{0.32,0.32,0.32}
\definecolor{gray33}{rgb}{0.33,0.33,0.33}
\definecolor{gray34}{rgb}{0.34,0.34,0.34}
\definecolor{gray35}{rgb}{0.35,0.35,0.35}
\definecolor{gray36}{rgb}{0.36,0.36,0.36}
\definecolor{gray37}{rgb}{0.37,0.37,0.37}
\definecolor{gray38}{rgb}{0.38,0.38,0.38}
\definecolor{gray39}{rgb}{0.39,0.39,0.39}
\definecolor{gray3}{rgb}{0.03,0.03,0.03}
\definecolor{gray40}{rgb}{0.40,0.40,0.40}
\definecolor{gray41}{rgb}{0.41,0.41,0.41}
\definecolor{gray42}{rgb}{0.42,0.42,0.42}
\definecolor{gray43}{rgb}{0.43,0.43,0.43}
\definecolor{gray44}{rgb}{0.44,0.44,0.44}
\definecolor{gray45}{rgb}{0.45,0.45,0.45}
\definecolor{gray46}{rgb}{0.46,0.46,0.46}
\definecolor{gray47}{rgb}{0.47,0.47,0.47}
\definecolor{gray48}{rgb}{0.48,0.48,0.48}
\definecolor{gray49}{rgb}{0.49,0.49,0.49}
\definecolor{gray4}{rgb}{0.04,0.04,0.04}
\definecolor{gray50}{rgb}{0.50,0.50,0.50}
\definecolor{gray51}{rgb}{0.51,0.51,0.51}
\definecolor{gray52}{rgb}{0.52,0.52,0.52}
\definecolor{gray53}{rgb}{0.53,0.53,0.53}
\definecolor{gray54}{rgb}{0.54,0.54,0.54}
\definecolor{gray55}{rgb}{0.55,0.55,0.55}
\definecolor{gray56}{rgb}{0.56,0.56,0.56}
\definecolor{gray57}{rgb}{0.57,0.57,0.57}
\definecolor{gray58}{rgb}{0.58,0.58,0.58}
\definecolor{gray59}{rgb}{0.59,0.59,0.59}
\definecolor{gray5}{rgb}{0.05,0.05,0.05}
\definecolor{gray60}{rgb}{0.60,0.60,0.60}
\definecolor{gray61}{rgb}{0.61,0.61,0.61}
\definecolor{gray62}{rgb}{0.62,0.62,0.62}
\definecolor{gray63}{rgb}{0.63,0.63,0.63}
\definecolor{gray64}{rgb}{0.64,0.64,0.64}
\definecolor{gray65}{rgb}{0.65,0.65,0.65}
\definecolor{gray66}{rgb}{0.66,0.66,0.66}
\definecolor{gray67}{rgb}{0.67,0.67,0.67}
\definecolor{gray68}{rgb}{0.68,0.68,0.68}
\definecolor{gray69}{rgb}{0.69,0.69,0.69}
\definecolor{gray6}{rgb}{0.06,0.06,0.06}
\definecolor{gray70}{rgb}{0.70,0.70,0.70}
\definecolor{gray71}{rgb}{0.71,0.71,0.71}
\definecolor{gray72}{rgb}{0.72,0.72,0.72}
\definecolor{gray73}{rgb}{0.73,0.73,0.73}
\definecolor{gray74}{rgb}{0.74,0.74,0.74}
\definecolor{gray75}{rgb}{0.75,0.75,0.75}
\definecolor{gray76}{rgb}{0.76,0.76,0.76}
\definecolor{gray77}{rgb}{0.77,0.77,0.77}
\definecolor{gray78}{rgb}{0.78,0.78,0.78}
\definecolor{gray79}{rgb}{0.79,0.79,0.79}
\definecolor{gray7}{rgb}{0.07,0.07,0.07}
\definecolor{gray80}{rgb}{0.80,0.80,0.80}
\definecolor{gray81}{rgb}{0.81,0.81,0.81}
\definecolor{gray82}{rgb}{0.82,0.82,0.82}
\definecolor{gray83}{rgb}{0.83,0.83,0.83}
\definecolor{gray84}{rgb}{0.84,0.84,0.84}
\definecolor{gray85}{rgb}{0.85,0.85,0.85}
\definecolor{gray86}{rgb}{0.86,0.86,0.86}
\definecolor{gray87}{rgb}{0.87,0.87,0.87}
\definecolor{gray88}{rgb}{0.88,0.88,0.88}
\definecolor{gray89}{rgb}{0.89,0.89,0.89}
\definecolor{gray8}{rgb}{0.08,0.08,0.08}
\definecolor{gray90}{rgb}{0.90,0.90,0.90}
\definecolor{gray91}{rgb}{0.91,0.91,0.91}
\definecolor{gray92}{rgb}{0.92,0.92,0.92}
\definecolor{gray93}{rgb}{0.93,0.93,0.93}
\definecolor{gray94}{rgb}{0.94,0.94,0.94}
\definecolor{gray95}{rgb}{0.95,0.95,0.95}
\definecolor{gray96}{rgb}{0.96,0.96,0.96}
\definecolor{gray97}{rgb}{0.97,0.97,0.97}
\definecolor{gray98}{rgb}{0.98,0.98,0.98}
\definecolor{gray99}{rgb}{0.99,0.99,0.99}
\definecolor{gray9}{rgb}{0.09,0.09,0.09}
\definecolor{gray}{rgb}{0.75,0.75,0.75}
\definecolor{green1}{rgb}{0.00,1.00,0.00}
\definecolor{green2}{rgb}{0.00,0.93,0.00}
\definecolor{green3}{rgb}{0.00,0.80,0.00}
\definecolor{green4}{rgb}{0.00,0.55,0.00}
\definecolor{greenyellow}{rgb}{0.68,1.00,0.18}
\definecolor{green}{rgb}{0.00,1.00,0.00}
\definecolor{grey0}{rgb}{0.00,0.00,0.00}
\definecolor{grey100}{rgb}{1.00,1.00,1.00}
\definecolor{grey10}{rgb}{0.10,0.10,0.10}
\definecolor{grey11}{rgb}{0.11,0.11,0.11}
\definecolor{grey12}{rgb}{0.12,0.12,0.12}
\definecolor{grey13}{rgb}{0.13,0.13,0.13}
\definecolor{grey14}{rgb}{0.14,0.14,0.14}
\definecolor{grey15}{rgb}{0.15,0.15,0.15}
\definecolor{grey16}{rgb}{0.16,0.16,0.16}
\definecolor{grey17}{rgb}{0.17,0.17,0.17}
\definecolor{grey18}{rgb}{0.18,0.18,0.18}
\definecolor{grey19}{rgb}{0.19,0.19,0.19}
\definecolor{grey1}{rgb}{0.01,0.01,0.01}
\definecolor{grey20}{rgb}{0.20,0.20,0.20}
\definecolor{grey21}{rgb}{0.21,0.21,0.21}
\definecolor{grey22}{rgb}{0.22,0.22,0.22}
\definecolor{grey23}{rgb}{0.23,0.23,0.23}
\definecolor{grey24}{rgb}{0.24,0.24,0.24}
\definecolor{grey25}{rgb}{0.25,0.25,0.25}
\definecolor{grey26}{rgb}{0.26,0.26,0.26}
\definecolor{grey27}{rgb}{0.27,0.27,0.27}
\definecolor{grey28}{rgb}{0.28,0.28,0.28}
\definecolor{grey29}{rgb}{0.29,0.29,0.29}
\definecolor{grey2}{rgb}{0.02,0.02,0.02}
\definecolor{grey30}{rgb}{0.30,0.30,0.30}
\definecolor{grey31}{rgb}{0.31,0.31,0.31}
\definecolor{grey32}{rgb}{0.32,0.32,0.32}
\definecolor{grey33}{rgb}{0.33,0.33,0.33}
\definecolor{grey34}{rgb}{0.34,0.34,0.34}
\definecolor{grey35}{rgb}{0.35,0.35,0.35}
\definecolor{grey36}{rgb}{0.36,0.36,0.36}
\definecolor{grey37}{rgb}{0.37,0.37,0.37}
\definecolor{grey38}{rgb}{0.38,0.38,0.38}
\definecolor{grey39}{rgb}{0.39,0.39,0.39}
\definecolor{grey3}{rgb}{0.03,0.03,0.03}
\definecolor{grey40}{rgb}{0.40,0.40,0.40}
\definecolor{grey41}{rgb}{0.41,0.41,0.41}
\definecolor{grey42}{rgb}{0.42,0.42,0.42}
\definecolor{grey43}{rgb}{0.43,0.43,0.43}
\definecolor{grey44}{rgb}{0.44,0.44,0.44}
\definecolor{grey45}{rgb}{0.45,0.45,0.45}
\definecolor{grey46}{rgb}{0.46,0.46,0.46}
\definecolor{grey47}{rgb}{0.47,0.47,0.47}
\definecolor{grey48}{rgb}{0.48,0.48,0.48}
\definecolor{grey49}{rgb}{0.49,0.49,0.49}
\definecolor{grey4}{rgb}{0.04,0.04,0.04}
\definecolor{grey50}{rgb}{0.50,0.50,0.50}
\definecolor{grey51}{rgb}{0.51,0.51,0.51}
\definecolor{grey52}{rgb}{0.52,0.52,0.52}
\definecolor{grey53}{rgb}{0.53,0.53,0.53}
\definecolor{grey54}{rgb}{0.54,0.54,0.54}
\definecolor{grey55}{rgb}{0.55,0.55,0.55}
\definecolor{grey56}{rgb}{0.56,0.56,0.56}
\definecolor{grey57}{rgb}{0.57,0.57,0.57}
\definecolor{grey58}{rgb}{0.58,0.58,0.58}
\definecolor{grey59}{rgb}{0.59,0.59,0.59}
\definecolor{grey5}{rgb}{0.05,0.05,0.05}
\definecolor{grey60}{rgb}{0.60,0.60,0.60}
\definecolor{grey61}{rgb}{0.61,0.61,0.61}
\definecolor{grey62}{rgb}{0.62,0.62,0.62}
\definecolor{grey63}{rgb}{0.63,0.63,0.63}
\definecolor{grey64}{rgb}{0.64,0.64,0.64}
\definecolor{grey65}{rgb}{0.65,0.65,0.65}
\definecolor{grey66}{rgb}{0.66,0.66,0.66}
\definecolor{grey67}{rgb}{0.67,0.67,0.67}
\definecolor{grey68}{rgb}{0.68,0.68,0.68}
\definecolor{grey69}{rgb}{0.69,0.69,0.69}
\definecolor{grey6}{rgb}{0.06,0.06,0.06}
\definecolor{grey70}{rgb}{0.70,0.70,0.70}
\definecolor{grey71}{rgb}{0.71,0.71,0.71}
\definecolor{grey72}{rgb}{0.72,0.72,0.72}
\definecolor{grey73}{rgb}{0.73,0.73,0.73}
\definecolor{grey74}{rgb}{0.74,0.74,0.74}
\definecolor{grey75}{rgb}{0.75,0.75,0.75}
\definecolor{grey76}{rgb}{0.76,0.76,0.76}
\definecolor{grey77}{rgb}{0.77,0.77,0.77}
\definecolor{grey78}{rgb}{0.78,0.78,0.78}
\definecolor{grey79}{rgb}{0.79,0.79,0.79}
\definecolor{grey7}{rgb}{0.07,0.07,0.07}
\definecolor{grey80}{rgb}{0.80,0.80,0.80}
\definecolor{grey81}{rgb}{0.81,0.81,0.81}
\definecolor{grey82}{rgb}{0.82,0.82,0.82}
\definecolor{grey83}{rgb}{0.83,0.83,0.83}
\definecolor{grey84}{rgb}{0.84,0.84,0.84}
\definecolor{grey85}{rgb}{0.85,0.85,0.85}
\definecolor{grey86}{rgb}{0.86,0.86,0.86}
\definecolor{grey87}{rgb}{0.87,0.87,0.87}
\definecolor{grey88}{rgb}{0.88,0.88,0.88}
\definecolor{grey89}{rgb}{0.89,0.89,0.89}
\definecolor{grey8}{rgb}{0.08,0.08,0.08}
\definecolor{grey90}{rgb}{0.90,0.90,0.90}
\definecolor{grey91}{rgb}{0.91,0.91,0.91}
\definecolor{grey92}{rgb}{0.92,0.92,0.92}
\definecolor{grey93}{rgb}{0.93,0.93,0.93}
\definecolor{grey94}{rgb}{0.94,0.94,0.94}
\definecolor{grey95}{rgb}{0.95,0.95,0.95}
\definecolor{grey96}{rgb}{0.96,0.96,0.96}
\definecolor{grey97}{rgb}{0.97,0.97,0.97}
\definecolor{grey98}{rgb}{0.98,0.98,0.98}
\definecolor{grey99}{rgb}{0.99,0.99,0.99}
\definecolor{grey9}{rgb}{0.09,0.09,0.09}
\definecolor{grey}{rgb}{0.75,0.75,0.75}
\definecolor{honeydew1}{rgb}{0.94,1.00,0.94}
\definecolor{honeydew2}{rgb}{0.88,0.93,0.88}
\definecolor{honeydew3}{rgb}{0.76,0.80,0.76}
\definecolor{honeydew4}{rgb}{0.51,0.55,0.51}
\definecolor{honeydew}{rgb}{0.94,1.00,0.94}
\definecolor{hotpink}{rgb}{1.00,0.41,0.71}
\definecolor{indianred}{rgb}{0.80,0.36,0.36}
\definecolor{ivory1}{rgb}{1.00,1.00,0.94}
\definecolor{ivory2}{rgb}{0.93,0.93,0.88}
\definecolor{ivory3}{rgb}{0.80,0.80,0.76}
\definecolor{ivory4}{rgb}{0.55,0.55,0.51}
\definecolor{ivory}{rgb}{1.00,1.00,0.94}
\definecolor{khaki1}{rgb}{1.00,0.96,0.56}
\definecolor{khaki2}{rgb}{0.93,0.90,0.52}
\definecolor{khaki3}{rgb}{0.80,0.78,0.45}
\definecolor{khaki4}{rgb}{0.55,0.53,0.31}
\definecolor{khaki}{rgb}{0.94,0.90,0.55}
\definecolor{lavenderblush}{rgb}{1.00,0.94,0.96}
\definecolor{lavender}{rgb}{0.90,0.90,0.98}
\definecolor{lawngreen}{rgb}{0.49,0.99,0.00}
\definecolor{lemonchiffon}{rgb}{1.00,0.98,0.80}
\definecolor{lightblue}{rgb}{0.68,0.85,0.90}
\definecolor{lightcoral}{rgb}{0.94,0.50,0.50}
\definecolor{lightcyan}{rgb}{0.88,1.00,1.00}
\definecolor{lightgoldenrod}{rgb}{0.93,0.87,0.51}
\definecolor{lightgoldenrod}{rgb}{0.98,0.98,0.82}
\definecolor{lightgray}{rgb}{0.83,0.83,0.83}
\definecolor{lightgreen}{rgb}{0.56,0.93,0.56}
\definecolor{lightgrey}{rgb}{0.83,0.83,0.83}
\definecolor{lightpink}{rgb}{1.00,0.71,0.76}
\definecolor{lightsalmon}{rgb}{1.00,0.63,0.48}
\definecolor{lightsea}{rgb}{0.13,0.70,0.67}
\definecolor{lightsky}{rgb}{0.53,0.81,0.98}
\definecolor{lightslate}{rgb}{0.47,0.53,0.60}
\definecolor{lightslate}{rgb}{0.47,0.53,0.60}
\definecolor{lightslate}{rgb}{0.52,0.44,1.00}
\definecolor{lightsteel}{rgb}{0.69,0.77,0.87}
\definecolor{lightyellow}{rgb}{1.00,1.00,0.88}
\definecolor{limegreen}{rgb}{0.20,0.80,0.20}
\definecolor{linen}{rgb}{0.98,0.94,0.90}
\definecolor{magenta1}{rgb}{1.00,0.00,1.00}
\definecolor{magenta2}{rgb}{0.93,0.00,0.93}
\definecolor{magenta3}{rgb}{0.80,0.00,0.80}
\definecolor{magenta4}{rgb}{0.55,0.00,0.55}
\definecolor{magenta}{rgb}{1.00,0.00,1.00}
\definecolor{maroon1}{rgb}{1.00,0.20,0.70}
\definecolor{maroon2}{rgb}{0.93,0.19,0.65}
\definecolor{maroon3}{rgb}{0.80,0.16,0.56}
\definecolor{maroon4}{rgb}{0.55,0.11,0.38}
\definecolor{maroon}{rgb}{0.69,0.19,0.38}
\definecolor{mediumaquamarine}{rgb}{0.40,0.80,0.67}
\definecolor{mediumblue}{rgb}{0.00,0.00,0.80}
\definecolor{mediumorchid}{rgb}{0.73,0.33,0.83}
\definecolor{mediumpurple}{rgb}{0.58,0.44,0.86}
\definecolor{mediumsea}{rgb}{0.24,0.70,0.44}
\definecolor{mediumslate}{rgb}{0.48,0.41,0.93}
\definecolor{mediumspring}{rgb}{0.00,0.98,0.60}
\definecolor{mediumturquoise}{rgb}{0.28,0.82,0.80}
\definecolor{mediumviolet}{rgb}{0.78,0.08,0.52}
\definecolor{midnightblue}{rgb}{0.10,0.10,0.44}
\definecolor{mintcream}{rgb}{0.96,1.00,0.98}
\definecolor{mistyrose}{rgb}{1.00,0.89,0.88}
\definecolor{moccasin}{rgb}{1.00,0.89,0.71}
\definecolor{navajowhite}{rgb}{1.00,0.87,0.68}
\definecolor{navyblue}{rgb}{0.00,0.00,0.50}
\definecolor{navy}{rgb}{0.00,0.00,0.50}
\definecolor{oldlace}{rgb}{0.99,0.96,0.90}
\definecolor{olivedrab}{rgb}{0.42,0.56,0.14}
\definecolor{orange1}{rgb}{1.00,0.65,0.00}
\definecolor{orange2}{rgb}{0.93,0.60,0.00}
\definecolor{orange3}{rgb}{0.80,0.52,0.00}
\definecolor{orange4}{rgb}{0.55,0.35,0.00}
\definecolor{orangered}{rgb}{1.00,0.27,0.00}
\definecolor{orange}{rgb}{1.00,0.65,0.00}
\definecolor{orchid1}{rgb}{1.00,0.51,0.98}
\definecolor{orchid2}{rgb}{0.93,0.48,0.91}
\definecolor{orchid3}{rgb}{0.80,0.41,0.79}
\definecolor{orchid4}{rgb}{0.55,0.28,0.54}
\definecolor{orchid}{rgb}{0.85,0.44,0.84}
\definecolor{palegoldenrod}{rgb}{0.93,0.91,0.67}
\definecolor{palegreen}{rgb}{0.60,0.98,0.60}
\definecolor{paleturquoise}{rgb}{0.69,0.93,0.93}
\definecolor{paleviolet}{rgb}{0.86,0.44,0.58}
\definecolor{papayawhip}{rgb}{1.00,0.94,0.84}
\definecolor{peachpuff}{rgb}{1.00,0.85,0.73}
\definecolor{peru}{rgb}{0.80,0.52,0.25}
\definecolor{pink1}{rgb}{1.00,0.71,0.77}
\definecolor{pink2}{rgb}{0.93,0.66,0.72}
\definecolor{pink3}{rgb}{0.80,0.57,0.62}
\definecolor{pink4}{rgb}{0.55,0.39,0.42}
\definecolor{pink}{rgb}{1.00,0.75,0.80}
\definecolor{plum1}{rgb}{1.00,0.73,1.00}
\definecolor{plum2}{rgb}{0.93,0.68,0.93}
\definecolor{plum3}{rgb}{0.80,0.59,0.80}
\definecolor{plum4}{rgb}{0.55,0.40,0.55}
\definecolor{plum}{rgb}{0.87,0.63,0.87}
\definecolor{powderblue}{rgb}{0.69,0.88,0.90}
\definecolor{purple1}{rgb}{0.61,0.19,1.00}
\definecolor{purple2}{rgb}{0.57,0.17,0.93}
\definecolor{purple3}{rgb}{0.49,0.15,0.80}
\definecolor{purple4}{rgb}{0.33,0.10,0.55}
\definecolor{purple}{rgb}{0.63,0.13,0.94}
\definecolor{red1}{rgb}{1.00,0.00,0.00}
\definecolor{red2}{rgb}{0.93,0.00,0.00}
\definecolor{red3}{rgb}{0.80,0.00,0.00}
\definecolor{red4}{rgb}{0.55,0.00,0.00}
\definecolor{red}{rgb}{1.00,0.00,0.00}
\definecolor{rosybrown}{rgb}{0.74,0.56,0.56}
\definecolor{royalblue}{rgb}{0.25,0.41,0.88}
\definecolor{saddlebrown}{rgb}{0.55,0.27,0.07}
\definecolor{salmon1}{rgb}{1.00,0.55,0.41}
\definecolor{salmon2}{rgb}{0.93,0.51,0.38}
\definecolor{salmon3}{rgb}{0.80,0.44,0.33}
\definecolor{salmon4}{rgb}{0.55,0.30,0.22}
\definecolor{salmon}{rgb}{0.98,0.50,0.45}
\definecolor{sandybrown}{rgb}{0.96,0.64,0.38}
\definecolor{seagreen}{rgb}{0.18,0.55,0.34}
\definecolor{seashell1}{rgb}{1.00,0.96,0.93}
\definecolor{seashell2}{rgb}{0.93,0.90,0.87}
\definecolor{seashell3}{rgb}{0.80,0.77,0.75}
\definecolor{seashell4}{rgb}{0.55,0.53,0.51}
\definecolor{seashell}{rgb}{1.00,0.96,0.93}
\definecolor{sienna1}{rgb}{1.00,0.51,0.28}
\definecolor{sienna2}{rgb}{0.93,0.47,0.26}
\definecolor{sienna3}{rgb}{0.80,0.41,0.22}
\definecolor{sienna4}{rgb}{0.55,0.28,0.15}
\definecolor{sienna}{rgb}{0.63,0.32,0.18}
\definecolor{skyblue}{rgb}{0.53,0.81,0.92}
\definecolor{slateblue}{rgb}{0.42,0.35,0.80}
\definecolor{slategray}{rgb}{0.44,0.50,0.56}
\definecolor{slategrey}{rgb}{0.44,0.50,0.56}
\definecolor{snow1}{rgb}{1.00,0.98,0.98}
\definecolor{snow2}{rgb}{0.93,0.91,0.91}
\definecolor{snow3}{rgb}{0.80,0.79,0.79}
\definecolor{snow4}{rgb}{0.55,0.54,0.54}
\definecolor{snow}{rgb}{1.00,0.98,0.98}
\definecolor{springgreen}{rgb}{0.00,1.00,0.50}
\definecolor{steelblue}{rgb}{0.27,0.51,0.71}
\definecolor{tan1}{rgb}{1.00,0.65,0.31}
\definecolor{tan2}{rgb}{0.93,0.60,0.29}
\definecolor{tan3}{rgb}{0.80,0.52,0.25}
\definecolor{tan4}{rgb}{0.55,0.35,0.17}
\definecolor{tan}{rgb}{0.82,0.71,0.55}
\definecolor{thistle1}{rgb}{1.00,0.88,1.00}
\definecolor{thistle2}{rgb}{0.93,0.82,0.93}
\definecolor{thistle3}{rgb}{0.80,0.71,0.80}
\definecolor{thistle4}{rgb}{0.55,0.48,0.55}
\definecolor{thistle}{rgb}{0.85,0.75,0.85}
\definecolor{tomato1}{rgb}{1.00,0.39,0.28}
\definecolor{tomato2}{rgb}{0.93,0.36,0.26}
\definecolor{tomato3}{rgb}{0.80,0.31,0.22}
\definecolor{tomato4}{rgb}{0.55,0.21,0.15}
\definecolor{tomato}{rgb}{1.00,0.39,0.28}
\definecolor{turquoise1}{rgb}{0.00,0.96,1.00}
\definecolor{turquoise2}{rgb}{0.00,0.90,0.93}
\definecolor{turquoise3}{rgb}{0.00,0.77,0.80}
\definecolor{turquoise4}{rgb}{0.00,0.53,0.55}
\definecolor{turquoise}{rgb}{0.25,0.88,0.82}
\definecolor{violetred}{rgb}{0.82,0.13,0.56}
\definecolor{violet}{rgb}{0.93,0.51,0.93}
\definecolor{wheat1}{rgb}{1.00,0.91,0.73}
\definecolor{wheat2}{rgb}{0.93,0.85,0.68}
\definecolor{wheat3}{rgb}{0.80,0.73,0.59}
\definecolor{wheat4}{rgb}{0.55,0.49,0.40}
\definecolor{wheat}{rgb}{0.96,0.87,0.70}
\definecolor{whitesmoke}{rgb}{0.96,0.96,0.96}
\definecolor{white}{rgb}{1.00,1.00,1.00}
\definecolor{yellow1}{rgb}{1.00,1.00,0.00}
\definecolor{yellow2}{rgb}{0.93,0.93,0.00}
\definecolor{yellow3}{rgb}{0.80,0.80,0.00}
\definecolor{yellow4}{rgb}{0.55,0.55,0.00}
\definecolor{yellowgreen}{rgb}{0.60,0.80,0.20}
\definecolor{yellow}{rgb}{1.00,1.00,0.00}
 
\usepackage[usenames,dvipsnames]{xcolor}
\usepackage[psamsfonts]{amsfonts}
\usepackage{graphicx}
\usepackage{amssymb}
\usepackage{shortvrb} 
 
\usepackage{hyperref} 
\hypersetup{pdfpagemode=FullScreen,   
backref,
colorlinks=true,  
citecolor=Bittersweet,
linkcolor=Bittersweet, 
urlcolor=Bittersweet
}
 
\setattribute{journal}{name}{}

\RequirePackage{hypernat}

\usepackage{hyperref} 
\hypersetup{pdfpagemode=FullScreen,   
backref,
colorlinks=true, 
citecolor=Bittersweet,
linkcolor=Bittersweet, 
urlcolor=Bittersweet
}

\usepackage[psamsfonts]{amsfonts}

\usepackage{pdflscape}

\usepackage{dsfont}
 
\defcitealias{Cut2015}{CPV15}

\newtheorem{Lem}{Lemma}[section]

\newtheorem{Theor}{Theorem}[section]
\newtheorem{Corol}{Corollary}[section]

\newcommand{\cqfd}{\hfill $\square$}

\newcommand{\R}{\mathbb R}

\newcommand{\n}{^{(n)}}

\newcommand{\Xb}{\mathbf{X}}
\newcommand{\Zb}{\mathbf{Z}}
\newcommand{\zb}{\mathbf{z}}
\newcommand{\Db}{\mathbf{D}}
\newcommand{\Jb}{\mathbf{J}}

\newcommand{\Wb}{\mathbf{W}}
\newcommand{\Bb}{\mathbf{B}}
\newcommand{\Ab}{\mathbf{A}}
\newcommand{\Sb}{\mathbf{S}}
\newcommand{\Sigb}{{\pmb \Sigma}}
\newcommand{\Lamb}{{\pmb \Lambda}}
\newcommand{\Gamb}{{\pmb \Gamma_p}}
\newcommand{\Deltab}{{\pmb \Delta}}

\newcommand{\xb}{\mathbf{x}}
\newcommand{\taub}{{\pmb\tau}}
\newcommand{\thetab}{{\pmb\theta}}
\newcommand{\varthetab}{{\pmb\vartheta}}
\newcommand{\pr}{^{\prime}}

\begin{document}
\begin{frontmatter}

\title{On the power of axial tests of uniformity on spheres}
\runtitle{On the power of axial tests of uniformity on spheres}

\begin{aug}
\author{\fnms{Christine} \snm{cutting}$^\dagger$, \fnms{Davy} \snm{Paindaveine}\thanksref{t1} and \fnms{Thomas} \snm{Verdebout}$^\ddagger$}

\thankstext{t1}{Corresponding author.}
\runauthor{Chr. Cutting, D. Paindaveine and Th. Verdebout}

\affiliation{Universit\'{e} libre de Bruxelles}

\address{$^{\dagger * \ddagger}$Universit\'{e} libre de Bruxelles\\
ECARES and D\'{e}partement de Math\'{e}matique\\
Avenue F.D. Roosevelt, 50\\
ECARES, CP114/04\\
B-1050, Brussels\\
Belgium\\
}

\address{$^*$Universit\'{e} Toulouse Capitole\\
Toulouse School of Economics\\
21, All\'{e}e de Brienne\\
31015 Toulouse Cedex 6\\
France\\
}

\end{aug}
\vspace{3mm}

\begin{abstract}
Testing uniformity on the $p$-dimensional unit sphere is arguably the most fundamental problem in directional statistics. In this paper, we consider this problem in the framework of \emph{axial} data, that is, under the assumption that the~$n$ observations at hand are randomly drawn from a distribution that charges antipodal regions equally. More precisely, we focus on axial, rotationally symmetric, alternatives and first address the problem under which the direction~$\thetab$ of the corresponding symmetry axis is specified. In this setup, we obtain Le Cam optimal tests of uniformity, that are based on the sample covariance matrix (unlike their non-axial analogs, that are based on the sample average). For the more important unspecified-$\thetab$ problem, some classical tests are available in the literature, but virtually nothing is known on their non-null behavior. We therefore study the non-null behavior of the celebrated Bingham test and of other tests that exploit the single-spiked nature of the considered alternatives. We perform Monte Carlo exercises to investigate the finite-sample behavior of our tests and to show their agreement with our asymptotic results.
%
\end{abstract}

\begin{keyword}[class=MSC]
\kwd[Primary ]{62H11, 62F05}
\kwd[; secondary ]{62E20}
\end{keyword}

\begin{keyword}
\kwd{Axial data}
\kwd{Contiguity} 
\kwd{Directional statistics}
\kwd{Local asymptotic normality} 
\kwd{Rotationally symmetric distributions}
\kwd{Tests of uniformity}
\end{keyword}

\end{frontmatter}


\section{Introduction} 
\label{sec:intro}

Directional statistics are concerned with data taking values on the unit hypersphere $\mathcal{S}^{p-1}:=\{\xb\in\R^p:\|\xb\|^2:=\xb'\xb=1\}$ of $\R^p$. Classical applications, that most often relate to the circular case ($p=2$) or spherical one ($p=3$), involve wind and animal migration data or belong to fields such as geology, paleomagnetism, or cosmology.  We refer, e.g., to \cite{Fish87}, \cite{MarJup2000}, and \cite{LV17book} for book-length treatments of the topic and for further applications. 

Arguably the most fundamental problem in directional statistics is the problem of testing for uniformity, which, for a random sample~$\Xb_{1},\ldots,\Xb_{n}$ at hand, consists in testing the null hypothesis that the observations are sampled from the uniform distribution over~$\mathcal{S}^{p-1}$. This is a very classical problem in multivariate analysis that can be traced back to \cite{Bernoulli1735}. As explained in the review paper \cite{GV19}, the topic has recently received a lot of attention: to cite only a few contributions, \cite{Ju08} proposed data-driven Sobolev tests, \cite{CuA09} and \cite{GNC19} proposed tests based on random projections, \cite{LPthanh14} studied the problem for noisy data, \cite{PaiVer2016} obtained the high-dimensional limiting behavior of some classical test statistics under the null hypothesis while \cite{GPV19} transformed some uniformity tests into tests of rotational symmetry. 

A classical uniformity test dates back to \cite{Ray1919} and rejects the null hypothesis for large values of~$\|\bar{\Xb}_n\|$, where~$\bar{\Xb}_n=\frac{1}{n}\sum_{i=1}^n \Xb_i$ is the sample mean of the observations. This test will detect only alternatives  whose mean vectors are non-zero, hence is therefore typically used when possible deviations from uniformity are suspected to be asymmetric. 
In particular, the Rayleigh test will show no power when the common distribution of the~$\Xb_i$'s is an antipodally symmetric distribution on the sphere, that is, when this distribution attributes the same probability to antipodal regions. Far from being the exception, such antipodally symmetric distributions are actually those that need be considered when practitioners are facing \emph{axial data}, that is, when one does not observe genuine locations on the sphere but rather axes (a typical example of axial data relates to the directions of optical axes in quartz crystals; see, e.g., \citealp{MarJup2000}). Models and inference for axial data have been considered a lot in the literature: to cite a few, \cite{Tyl87ang} and more recently \cite{PRV19} considered inference on the distribution of the spatial sign of a Gaussian vector, \cite{Watson65}, \cite{Bi07} and \cite{Sra13} considered inference for Watson distributions (see the next section), \cite{Dryden05} obtained distributions on high-dimensional spheres while  \cite{AndSte1972}, \cite{Bin1974} and \cite{JU01} considered uniformity tests against antipodally symmetric alternatives.

Now, while the literature offers both axial and non-axial tests of uniformity, axial procedures unfortunately remain much less well understood than their non-axial counterparts, particularly so when it comes to their non-null behaviors. Strong results have been obtained for non-axial tests of uniformity regarding their asymptotic power under suitable local alternatives to uniformity and even regarding their optimality (we refer to \citealp{PaiVer17a} and to the references therein), but virtually nothing is known in that direction for axial tests of uniformity. This provides the main motivation for the present work, that intends to fill an important gap by studying the non-null behavior of some classical (and less classical) axial tests of uniformity. Quite naturally, we will do so in the semiparametric distributional framework that has been classically considered for non-axial tests of uniformity, namely the framework of \emph{rotationally symmetric distributions} indexed by a finite-dimensional parameter $(\kappa, \thetab) \in \R^+ \times {\cal S}^{p-1}$ and an infinite-dimensional parameter $f \in {\cal F}$ (a family of functions we  define in the next section). Within this semiparametric model, the null hypothesis of uniformity takes the form ${\cal H}_0: \kappa=0$. In this paper, we first derive the shape of tests that are locally and asymptotically optimal under any $f \in {\cal F}$ for the specified-$\thetab$ problem. We then focus on the unspecified-$\thetab$ problem and discuss its connection with the specified-$\thetab$ problem. We derive the limiting behavior, under sequences of contiguous alternatives (for any $f \in {\cal F}$), of the tests provided in  \cite{AndSte1972} and \cite{Bin1974}. Doing so, we obtain in particular the limiting behavior, under local alternatives, of the extreme eigenvalues of the spatial sign covariance matrix, which is a result of independent interest; see \cite{DuTyl16} and the references therein for a recent study of these eigenvalues.

The outline and contribution of the paper are as follows. 
In Section~\ref{sec:axial}, we introduce the class of (rotationally symmetric) alternatives to uniformity that will be considered in this work, which deviate from uniformity along a direction~$\thetab$ and have a severity controlled by a concentration parameter~$\kappa$. In this framework, we identify the sequences~$(\kappa_n)$ that make the corresponding sequences of alternatives contiguous to the null hypothesis of uniformity. 
In Section~\ref{sec:TestSpeci}, we tackle the problem of testing uniformity under specified~$\thetab$ and show that the resulting model is \emph{locally asymptotically normal} (LAN). We define the resulting optimal tests of uniformity and determine their asymptotic powers under contiguous alternatives.  
In Section~\ref{sec:Bingham}, we turn to the unspecified-$\thetab$ problem, we show that our LAN result naturally leads to the  \cite{Bin1974} test of uniformity, and we study the limiting behavior of this test under contiguous alternatives. 
In Section~\ref{sec:eigentest}, we turn our attention to tests that take into account the ``single-spiked" structure of the considered alternatives.
We characterize the asymptotic behavior of these tests both under the null hypothesis and under sequences of contiguous alternatives.   
While all results above are confirmed by suitable numerical exercises in Sections~\ref{sec:TestSpeci}--\ref{sec:eigentest}, we specifically conduct, in Section~\ref{sec:Simu}, Monte-Carlo simulations in order to compare the finite-sample powers of the various tests.  
%
We provide final comments and perspectives for future research in Section~\ref{sec:Final}. 
Finally, an appendix contains all proofs.   
%


\section{Axial rotationally symmetric distributions}
\label{sec:axial}

In this section, we set the notation and describe the class of axial distributions we will use as alternatives to uniformity. A celebrated class of axial distributions on the sphere is the one that collects the Watson distributions, which admit a density (throughout, densities over~$\mathcal{S}^{p-1}$ are with respect to the surface area measure) of the form 
\begin{equation}
\label{densWatson}
\xb\mapsto 
\frac{c_{p,\kappa} \Gamma(\frac{p-1}{2})}{2\pi^{(p-1)/2}} 
\exp(\kappa\, (\xb'\thetab)^2)
,
\qquad
\xb\in\mathcal{S}^{p-1}
,
\end{equation}
where~$\thetab$ belongs to~$\mathcal{S}^{p-1}$, $\kappa$ is a real number, and where the value of the normalizing constant~$c_{p,\kappa}$ can be obtained from~(\ref{normconst}); as usual, $\Gamma$ denotes Euler's Gamma function. Since this density is a symmetric function of~$\xb$, it attributes the same probability to antipodal regions on the sphere, hence is indeed suitable for axial data. The Watson distribution is rotationally symmetric about~$\thetab$, in the sense that if~$\Xb$ has density~(\ref{densWatson}), then ${\bf O}\Xb$ and~$\Xb$ share the same distribution for any $p\times p$ orthogonal matrix~${\bf O}$ such that~${\bf O}\thetab=\thetab$; consequently, $\thetab$ will be considered a \emph{location} parameter. The Watson distributions are the rotationally symmetric (or single-spiked) \cite{Bin1974} distributions. The parameter~$\kappa$ is a \emph{concentration} parameter: the larger~$|\kappa|$, the more the probability mass will be  concentrated---symmetrically about the poles~$\pm\thetab$ for positive values of~$\kappa$ (bipolar case) or symmetrically about the hyperspherical equator~$\mathcal{S}^\perp_{\thetab}:=\{\xb\in\mathcal{S}^{p-1}:\xb'\thetab=0\}$ for negative values of~$\kappa$ (girdle case). Of course, the value~$\kappa=0$ corresponds to the uniform distribution over the sphere. 

In this paper, we consider a natural semiparametric extension of the class of Watson distributions, namely the class of axial distributions admitting a density of the form 
\begin{equation}
\label{densconc}
\xb\mapsto 
\frac{c_{p,\kappa,f} \Gamma(\frac{p-1}{2})}{2\pi^{(p-1)/2}} 
f(\kappa\, (\xb'\thetab)^2)
,
\qquad
\xb\in\mathcal{S}^{p-1},
\end{equation}
where~$\thetab$ and~$\kappa$ are as in~(\ref{densWatson}), $f$ belongs to the class of functions~$\mathcal{F}:=\{f:\R\to\R^+: f \textrm{ monotone increasing, twice differentiable at }0, \textrm{ with } f(0)=f'(0)=1\}$, and where
\begin{equation}
	\label{normconst}
	c_{p,\kappa,f}
	=
	1\ \Big/ \int_{-1}^1 (1-s^2)^{(p-3)/2}f(\kappa s^2)\,ds
	. 
\end{equation}
The parameter~$\kappa$ is still a concentration parameter, that shares the same interpretation as for Watson distributions. The restriction to~$\mathcal{F}$ above is made for identifiability purposes. If~$\kappa\neq 0$, then~$f$ and the pair~$\{\pm\thetab\}$ are identifiable, but~$\thetab$ itself is not (which is natural for axial distributions). For~$\kappa=0$, the uniform distribution over the sphere is obtained (irrespective of~$f$), in which case the location parameter~$\thetab$ is unidentifiable, even up to a sign. The corresponding normalizing constant is~$c_p :=\lim_{\kappa\to 0} c_{p,\kappa,f}$. The distribution associated with~(\ref{densconc}) is rotationally symmetric about~$\thetab$, and if~$\Xb$ is a random vector with this distribution, then~$\Xb'\thetab$ has density~$s\mapsto c_{p,\kappa,f}(1-s^2)^{(p-3)/2}f(\kappa s^2)\mathbb{I}[|s|\leq 1]$, which explains the expression~(\ref{normconst}). In the present axial case, this density is of course symmetric with respect to zero. 

The semiparametric class of distributions just introduced will be used in the paper as alternatives to uniformity on the sphere. We will consider the following hypotheses and asymptotic scenarios. 
\vspace{-.5mm}
For~$\thetab\in\mathcal{S}^{p-1}$, a sequence~$(\kappa_n)$ in~$\R_0$ and~$f\in\mathcal{F}$, we will denote as~${\rm P}\n_{\thetab,\kappa_n,f}$ the hypothesis under which~$\Xb_{n1},\ldots,\Xb_{nn}$ form a random sample from the density~$\xb\mapsto c_{p,\kappa_n,f} f(\kappa_n\, (\xb'\thetab)^2)$ over~$\mathcal{S}^{p-1}$. This triangular array framework will allow us to consider local alternatives, associated with 
\vspace{-.5mm}
suitable sequences~$(\kappa_n)$ converging to zero. The null hypothesis of 
\vspace{-.5mm}
uniformity will be denoted as~${\rm P}\n_{0}$. The sequence of hypotheses~${\rm P}\n_{\thetab,\kappa_n,f}$ determines a sequence of alternatives to uniformity: the larger~$|\kappa_n|$, the more severe the corresponding alternative, whereas the sign of~$\kappa_n$ determines the type of alternatives considered, i.e., \emph{bipolar} (for~$\kappa_n>0$) or \emph{girdle-like} (for~$\kappa_n<0$). At places, it will be of interest to compare our results with those  obtained in the non-axial case, that is, in the case where~$\Xb_{n1},\ldots,\Xb_{nn}$ have a common density proportional to~$f(\kappa_n \xb'\thetab)$ still with~$f$ monotone increasing (rather than~$f(\kappa_n(\xb'\thetab)^2)$). 

Our first result identifies
\vspace{-1.1mm}
 the sequences of alternatives~${\rm P}\n_{\thetab,\kappa_n,f}$ that are contiguous to the sequence of null hypotheses~${\rm P}\n_{0}$.

\begin{Theor}
\label{TheorContig}
Fix~$p\in\{2,3,\ldots\}$, $\thetab\in\mathcal{S}^{p-1}$, and~$f\in\mathcal{F}$. Let~$(\kappa_n)$ be a sequence in~$\R_0$ that is~$O(1/\sqrt{n})$. 
\vspace{-.8mm}
Then, the sequence of alternative hypotheses~${\rm P}\n_{\thetab,\kappa_n,f}$ and the sequence of null hypotheses~${\rm P}\n_{0}$ are mutually contiguous.
\end{Theor}
 
In other words, if $\kappa_n=O(1/\sqrt{n})$, then no test for~$\mathcal{H}_{0n}:\{{\rm P}\n_{0}\}$ against~$\mathcal{H}_{1n}:\{{\rm P}\n_{\thetab,\kappa_n,f}\}$ can be consistent. Actually, it will follow from Theorem~\ref{TheorLAN} in the next section that~$1/\sqrt{n}$ is the \emph{contiguity rate}, in the sense that if~$\kappa_n=\tau/\sqrt{n}$ 
\vspace{-.3mm}
(for some non-zero real constant~$\tau$), then there exist tests for~$\mathcal{H}_{0n}:\{{\rm P}\n_{0}\}$ against~$\mathcal{H}_{1n}:\{{\rm P}\n_{\thetab,\kappa_n,f}\}$ 
\vspace{-.2mm}
showing non-trivial asymptotic powers (that is, asymptotic powers in~$(\alpha,1)$, where~$\alpha$ denotes the nominal level). The  contiguity rate in the axial case thus coincides with the one obtained in the non-axial case; see Theorems~2.1 and~3.1 in \cite{PaiVer17a}. 


\section{Tests of uniformity under specified location}
\label{sec:TestSpeci}

In this section, we consider the problem of testing uniformity over~$\mathcal{S}^{p-1}$ against the class of alternatives introduced in the previous section, in a situation where the location~$\thetab$ is specified. In other words, this corresponds to cases where it is known in which direction the possible deviation from uniformity would materialize. Depending on the exact type of alternatives we want to focus on (bipolar, girdle-type, or both), 
\vspace{-.5mm}
we will then consider, for a fixed~$\thetab$, the problem of testing~$\mathcal{H}_0\n:\{{\rm P}_0\n\}$ 
\vspace{-.8mm}
against (i) $\mathcal{H}_1\n:\cup_{\kappa>0}\cup_{f\in\mathcal{F}} \{{\rm P}_{\thetab,\kappa,f}\n\}$, against (ii) $\mathcal{H}_1\n:\cup_{\kappa<0}\cup_{f\in\mathcal{F}} \{{\rm P}_{\thetab,\kappa,f}\n\}$, or against (iii) $\mathcal{H}_1\n:\cup_{\kappa\neq 0}\cup_{f\in\mathcal{F}} \{{\rm P}_{\thetab,\kappa,f}\n\}$. 
%
Optimal testing may be based on the following Local Asymptotic Normality (LAN) result.

\begin{Theor}
\label{TheorLAN}
Fix~$p\in\{2,3,\ldots\}$, $\thetab\in\mathcal{S}^{p-1}$, and~$f\in\mathcal{F}$. Let~$\kappa_n=\tau_n p/\sqrt{n}$, where the real sequence~$(\tau_n)$ is~$O(1)$ but not~$o(1)$. Then,  
letting
$$
\Delta_{\thetab}\n
:=
\frac{p}{\sqrt{n}}
\sum_{i=1}^n \Big\{
(\Xb_{ni}\pr \thetab)^2-\frac{1}{p}
\Big\}
\quad\textrm{and}\quad
\Gamma_p 
:=
\frac{2(p-1)}{p+2}
\,
,
$$
we have that, as~$n\to\infty$ under~${\rm P}_{0}\n$,
\begin{equation}
\label{LAN}
\Lambda_n
=
\log \frac{d{\rm P}\n_{\thetab,\kappa_n,f}}{d{\rm P}\n_{0}} 
=
\tau_n
\Delta_{\thetab}\n
-
\frac{\tau_n^2}{2} \Gamma_p
+
o_{\rm P}(1)
,
\end{equation}
where~$\Delta_{\thetab}\n$
\vspace{-.5mm}
 is asymptotically normal with mean zero and variance~$\Gamma_p$. 
In other words, the sequence~$(\{ {\rm P}\n_{\thetab,\kappa,f} : \kappa\in\R \})$ 
\vspace{-.6mm}
$($where we let~${\rm P}\n_{\thetab,0,f}:={\rm P}\n_0)$ is 
\vspace{-.3mm}
locally asymptotically normal at~$\kappa=0$ with central sequence~$\Delta_{\thetab}\n$, Fisher information~$\Gamma_p$, and contiguity rate~$1/\sqrt{n}$.
\end{Theor}

This result confirms that~$1/\sqrt{n}$ is the contiguity rate when testing uniformity against the considered axial alternatives.
Note also that the central sequence rewrites
$$
\Delta_{\thetab}\n
=
\sqrt{n}
\big(
p
\,
\thetab'
\Sb_n
\thetab
-
1
\big)
,
$$ 
where~$\Sb_n:=n^{-1}\sum_{i=1}^n \Xb_{ni}\Xb_{ni}\pr$ is the sample covariance matrix of the observations (with respect to a fixed location, namely the origin of~$\R^{p}$). Consequently, optimal testing of uniformity for axial data will be based on~$\mathbf{S}_n$. This is to be compared with the non-axial case considered in \cite{PaiVer17a}, where optimal testing of uniformity is rather based on~$\bar{\Xb}_n=n^{-1}\sum_{i=1}^n \Xb_{ni}$. This will have important consequences when considering the unspecified-$\thetab$ case we turn to in Sections~\ref{sec:Bingham}--\ref{sec:eigentest}.

More importantly, the optimal axial tests of uniformity in the specified location case directly result from the LAN property above. 
\vspace{-.4mm}
More precisely, Theorem~\ref{TheorLAN} entails that, for the problem of testing~$\mathcal{H}_0\n:\{{\rm P}_0\n\}$ against $\cup_{\kappa>0}\cup_{f\in\mathcal{F}} \{{\rm P}_{\thetab,\kappa,f}\n\}$, the test~$\phi_{\thetab+}^{(n)}$ rejecting the null hypothesis at asymptotic level~$\alpha$ whenever 
\begin{equation}
	\label{Testright}
T\n_{\thetab}
:=
\frac{\Delta_{\thetab}\n}{\sqrt{\Gamma_p}}
>
z_{\alpha}
\end{equation}
is locally asymptotically most powerful; here, $z_\alpha=\Phi^{-1}(1-\alpha)$ denotes the upper \mbox{$\alpha$-quantile} of the standard normal distribution. 
\vspace{-.5mm}
A routine application of Le Cam's third lemma shows that, 
\vspace{-1.3mm}
under ${\rm P}\n_{\thetab,\kappa_n,f}$ with $\kappa_n=\tau p/\sqrt{n}$ ($\tau>0$),
$T\n_{\thetab}$ is asymptotically normal with mean $\Gamma_p^{1/2}\tau$ and variance one. 
\vspace{-.5mm}
Therefore, the corresponding asymptotic power of~$\phi_{\thetab+}^{(n)}$ is
\begin{equation}
\label{eq:right-sided-power}
\lim_{n\to\infty}
{\rm P}\n_{\thetab,\kappa_n,f}[T\n_{\thetab}>z_{\alpha}]
=
1-\Phi\big(z_{\alpha}-\Gamma_p^{1/2}\tau\big)
.
\end{equation}
Note that this asymptotic power does not converge to~$\alpha$ as~$p$ diverges to infinity. This may be surprising at first since departures from uniformity here are of a \emph{single-spiked} nature, that is, only materialize in a single direction out of the $p$ directions in~$\mathcal{S}^{p-1}$. The fact that this asymptotic power does not fade out for larger dimensions is actually explained by the fact that we did not consider local alternatives associated with~$\kappa_n=\tau/\sqrt{n}$ but rather with~$\kappa_n=\tau p/\sqrt{n}$, which properly scales local alternatives for different dimensions~$p$.

Optimal tests for the other 
\vspace{-.5mm}
one-sided problem and for the two-sided problem are obtained in a similar way. 
\vspace{-.5mm}
More precisely, for the problem of testing~$\mathcal{H}_0\n:\{{\rm P}_0\n\}$ against $\cup_{\kappa<0}\cup_{f\in\mathcal{F}} \{{\rm P}_{\thetab,\kappa,f}\n\}$, the test~$\phi_{\thetab-}^{(n)}$ rejecting 
\vspace{-.9mm}
the null hypothesis of uniformity at asymptotic level~$\alpha$ whenever~$T\n_{\thetab}<-z_\alpha$ is locally asymptotically most powerful and has asymptotic power
$$
\lim_{n\to\infty}
{\rm P}\n_{\thetab,\kappa_n,f}[T\n_{\thetab}< -z_{\alpha}]
=
\Phi(-z_{\alpha}-\Gamma_p^{1/2}\tau)
$$
under ${\rm P}\n_{\thetab,\kappa_n,f}$ with $\kappa_n=\tau p/\sqrt{n}$ ($\tau<0$). 
\vspace{-1.1mm}
The corresponding two-sided test, $\phi_{\thetab\pm}^{(n)}$ say, rejects the null hypothesis at asymptotic level~$\alpha$ whenever~$|T\n_{\thetab}|>z_{\alpha/2}$. This test is locally asymptotically maximin for the two-sided problem and has asymptotic power  
$$
\lim_{n\to\infty}{\rm P}\n_{\thetab,\kappa_n,f}[|T\n_{\thetab}|>z_{\alpha/2}]
=
2-\Phi(z_{\alpha/2}-\Gamma_p^{1/2}\tau)-\Phi(z_{\alpha/2}+\Gamma_p^{1/2}\tau)
$$
under ${\rm P}\n_{\thetab,\kappa_n,f}$ with $\kappa_n=\tau p/\sqrt{n}$ ($\tau\neq 0$). Again, the local asymptotic powers of these tests do not fade out for larger dimensions~$p$ but rather converge to a constant larger than~$\alpha$.

\begin{figure}
\centering
\includegraphics[width=\textwidth]{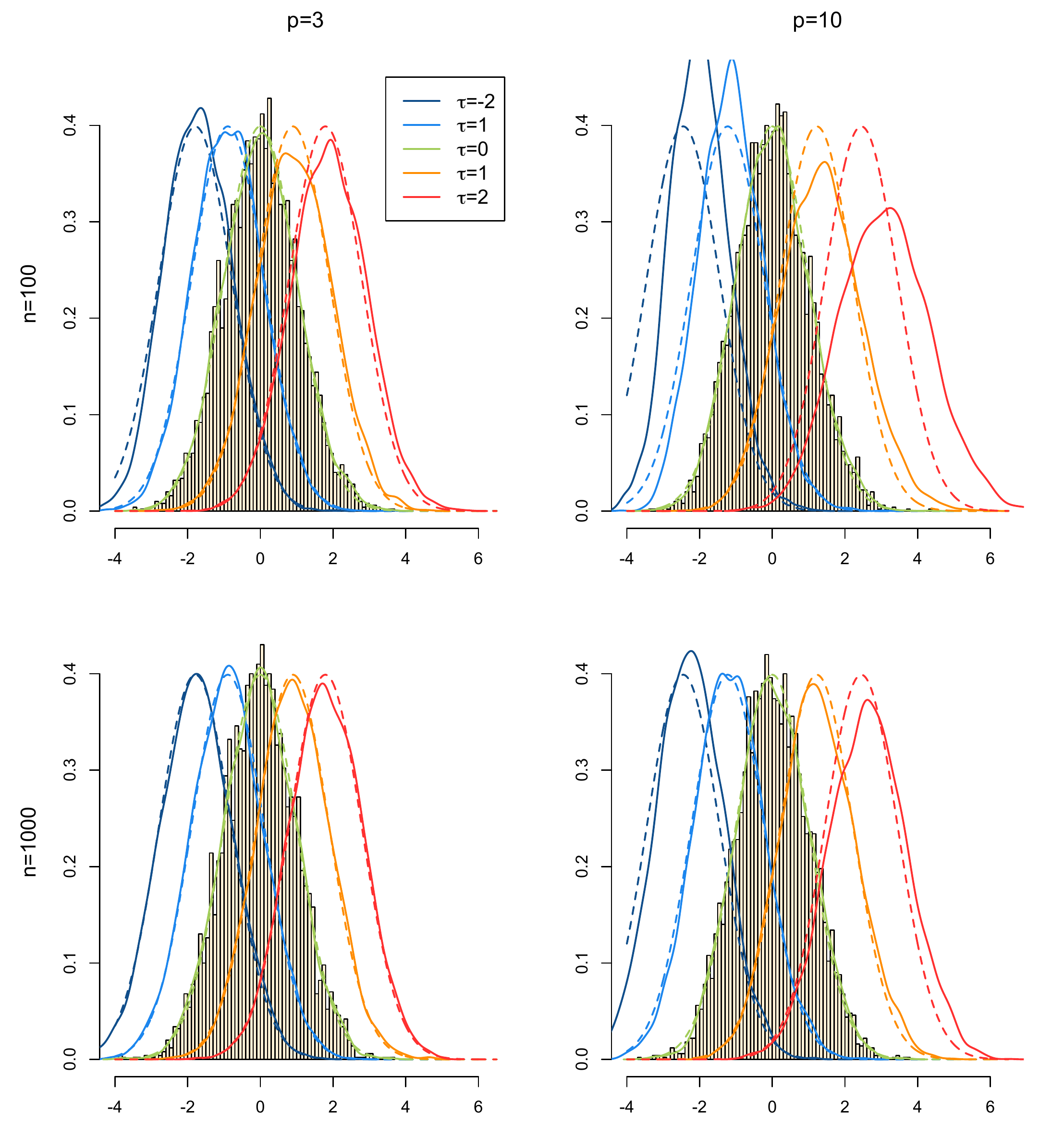}
\caption{
Plots of the kernel density estimates (solid curves) of the values of the test statistic~$T\n_{\thetab}$ in~(\ref{Testright}) obtained from~$M=5\:\!000$ independent random samples, of size~$n=100$ (top) or $1\:\!000$ (bottom), from the Watson distribution 
\vspace{-.6mm}
with location~$\thetab=(1,0,\ldots,0)'\in\R^p$ and concentration~$\kappa=\tau p/\sqrt{n}$, with~$\tau=-2,-1,0,1,2$ and with~$p=3$ (left) or~$p=10$ (right); 
\vspace{-.8mm}
 for~$\tau=0$, histograms of the values of~$T\n_{\thetab}$ are shown. The densities of the corresponding asymptotic ($\mathcal{N}(\Gamma_p^{1/2}\tau,1)$) distributions are also plotted (dashed curves). Throughout this paper, kernel density estimates are obtained from the R command \texttt{density} with default parameter values.
}
\label{Fig1}
\end{figure}

We conducted the following Monte Carlo exercise in order to check the validity of our asymptotic results. For any combination~$(n,p)$ of sample size~$n\in\{100,1\:\!000\}$ and dimension~$p\in\{3,10\}$, we generated collections of $5\:\!000$ independent random samples of size~$n$ from the Watson distribution with location~$\thetab=(1,0,\ldots,0)'\in\R^{p}$ and concentration~$\kappa_n=\tau p/\sqrt{n}$, for $\tau=-2,-1,0,1,2$; see Section~\ref{sec:axial}. The value~$\tau=0$ corresponds to the null hypothesis of uniformity over~$\mathcal{S}^{p-1}$, whereas the larger the non-zero value of~$|\tau|$ is, the more severe the alternative is. 
\vspace{-.5mm}
Kernel density estimates of the resulting values of the test statistic~$T\n_{\thetab}$ in~(\ref{Testright}) are provided in Figure \ref{Fig1}, that further plots the densities of the corresponding asymptotic distributions 
\vspace{-.5mm}
(for the null case~$\tau=0$, histograms of the values of~$T\n_{\thetab}$ are also shown). Clearly, our asymptotic results are confirmed by these simulations (yet, unsurprisingly, larger dimensions require larger sample sizes for asymptotic results to materialize). 

%


\section{The unspecified location case: the Bingham test}
\label{sec:Bingham}

We 
\vspace{-.4mm}
now turn to the \mbox{unspecified-$\thetab$} version of the testing problems considered in the previous section. 
\vspace{-.5mm}
We focus first on the one-sided problem of testing~$\mathcal{H}_0\n:\{{\rm P}_0\n\}$ against $\mathcal{H}_1\n:\cup_{\thetab\in\mathcal{S}^{p-1}}\cup_{\kappa>0}\cup_{f\in\mathcal{F}} \{{\rm P}_{\thetab,\kappa,f}\n\}$. 
It is 
\vspace{-.4mm}
convenient to reparameterize the submodel associated with~$\kappa\geq 0$ by writing~$\varthetab=\sqrt{\kappa} \thetab$. In this 
\vspace{-.4mm}
new parametrization (which, unlike the original
\vspace{-.5mm}
 curved one, is flat), the testing problem writes~$\mathcal{H}_0\n:\{{\rm P}_{\bf 0}\n={\rm P}_0\n\}$ against $\cup_{\varthetab\in\R^{p}\setminus\{0\}}\cup_{f\in\mathcal{F}} \{{\rm P}_{\varthetab,f}\n\}$, that is, 
\vspace{-.9mm}
simply consists in testing~$\mathcal{H}_0\n:\varthetab={\bf 0}$ against $\mathcal{H}_1\n:\varthetab\neq {\bf 0}$. Theorem~\ref{TheorLANunspec} below then describes the asymptotic behavior of the corresponding local log-likelihood ratios. 

To be able to state the result, we need to introduce the following notation. For a matrix~$\Ab$, we will write~${\rm vec}\, \Ab$ for the vector obtained by stacking the columns of~$\Ab$ on top of each other. We let~$\Jb_p:=({\rm vec}\, \mathbf{I}_p)({\rm vec}\, \mathbf{I}_p)'$, where~$\mathbf{I}_\ell$ is the $\ell\times \ell$ identity matrix. Finally, with the usual Kronecker product~$\otimes$, the $p^2\times p^2$ commutation matrix is~${\bf K}_p:=\sum_{i,j=1}^p ({\bf e}_i {\bf e}_j')\otimes ({\bf e}_j {\bf e}_i')$, where~$e_\ell$ is the $\ell$th vector of the canonical basis of~$\R^p$. We then have the following result.

\begin{Theor}
\label{TheorLANunspec}
Fix~$p\in\{2,3,\ldots\}$ and~$f\in\mathcal{F}$. Let~$\varthetab_n=(p/\sqrt{n})^{1/2}\taub_n$, where $(\taub_n)$ is a sequence in~$\R^p$ that is $O(1)$ but not~$o(1)$. Then, letting
$$
\Deltab\n
:=
p\sqrt{n}
\,
{\rm vec}\bigg(
\Sb_n-\frac{1}{p} \mathbf{I}_p
\bigg)
\quad\textrm{and}\quad
\Gamb
:=
\frac{p}{p+2} \Big( \mathbf{I}_{p^2} + \mathbf{K}_p -\frac{2}{p} \mathbf{J}_p \Big)
\,
,
$$
we have that, as~$n\to\infty$ under~${\rm P}_{0}\n$,
\begin{equation}
\label{LANvec}
\log \frac{d{\rm P}\n_{\varthetab_n,f}}{d{\rm P}\n_{0}} 
=
({\rm vec}(\taub_n\taub_n'))'
\Deltab\n
-
\frac{1}{2} ({\rm vec}(\taub_n\taub_n'))'\Gamb {\rm vec}(\taub_n\taub_n')
+
o_{\rm P}(1)
,
\end{equation}
where~$\Deltab\n$ is, still under~${\rm P}_{0}\n$, asymptotically normal with mean vector zero and covariance matrix~$\Gamb$.
\end{Theor}

Theorem~\ref{TheorLANunspec} shows that the contiguity rate for~$\varthetab$ is~$n^{-1/4}$, which corresponds to the contiguity rate~$n^{-1/2}$ obtained for~$\kappa$ in Theorem~\ref{TheorLAN} (recall that~$\varthetab=\sqrt{\kappa} \thetab$); however, as we will explain below, the limiting experiment in Theorem~\ref{TheorLANunspec} is non-standard. 
Writing~$\Ab^-$ for the Moore-Penrose generalized inverse of~$\Ab$, a natural test of uniformity is the test rejecting the null hypothesis at asymptotic level~$\alpha$ whenever
\begin{equation}
	\label{Binghamtest}
Q\n
:=
(\Deltab\n)' {\pmb\Gamma}_p^- \Deltab\n
=
\frac{np(p+2)}{2}
\bigg(\text{tr}[\Sb_n^2]-\frac{1}{p}\bigg)
>
\chi^2_{d_p,1-\alpha}
,
\end{equation}
where we denoted as~$\chi^2_{d_p,1-\alpha}$ the upper $\alpha$-quantile of the chi-square distribution with~$d_p:=p(p+1)/2-1$ degrees of freedom. 
\vspace{-.5mm}
This test, which rejects the null hypothesis when the sample variance of the eigenvalues~$\hat{\lambda}_{n1},\ldots,\hat{\lambda}_{np}$ of~$\Sb_n$ is too large, also addresses the problem of testing uniformity against the one-sided alternatives associated with~$\kappa<0$ or against the two-sided alternatives associated with~$\kappa\neq 0$. This procedure, which is known as the \cite{Bin1974} test (hence will be denoted as~$\phi_{\rm Bing}$ in the sequel), is often regarded as the simplest test of uniformity for axial data; see Section~10.7 in \cite{MarJup2000}. When applied to the unit vectors~$\Xb_{ni}=\Zb_{ni}/\|\Zb_{ni}\|$, $i=1,\ldots,n$, obtained from Euclidean data~$\Zb_{ni}$, $i=1,\ldots,n$, this test is also the sign test of sphericity from \cite{HalPai2006}, and it follows from that paper that, as an axial test for uniformity on the sphere, the Bingham test is optimal against \emph{angular Gaussian alternatives} (see \citealp{Tyl87ang}), that is, against projections of elliptical distributions on the sphere. 

Local asymptotic powers of the Bingham test can be obtained from the LAN result in Theorem~\ref{TheorLAN} and Le Cam's third lemma. We have the following result.

\begin{Theor}
\label{TheorBinghamHD}
Fix~$p\in\{2,3,\ldots\}$, $\thetab\in\mathcal{S}^{p-1}$, and~$f\in\mathcal{F}$. Let~$\kappa_n=\tau_n p/\sqrt{n}$, where the real sequence~$(\tau_n)$ converges to~$\tau$. Then, under ${\rm P}\n_{\thetab,\kappa_n,f}$,
\begin{equation}
	\label{Binghamtestncp}
Q\n
\stackrel{\mathcal{D}}{\to}
\chi^2_{d_p}
\bigg(
\frac{2(p-1)\tau^2}{p+2}
\bigg)
,
\end{equation}
where~$\chi^2_\ell(\delta)$ denotes the non-central chi-square distribution with~$\ell$ degrees of freedom and non-centrality parameter~$\delta$. Under the same sequence of alternatives, the asymptotic power of the Bingham test is therefore 
\begin{equation}
\label{eq:right-sided-powerBingham}
\lim_{n\to\infty}
{\rm P}\n_{\thetab,\kappa_n,f}[
Q\n
>
\chi^2_{d_p,1-\alpha}]
=
1
-
\Psi_{d_p}\bigg(
\chi^2_{d_p,1-\alpha}
;
\frac{2(p-1)\tau^2}{p+2}
\bigg)
,
\end{equation}
where~$\Psi_{\ell}(\cdot;\delta)$ is the cumulative distribution function of the~$\chi^2_\ell(\delta)$ distribution.
\end{Theor}

This result in particular shows that the Bingham test is a two-sided procedure, as the asymptotic power in~(\ref{eq:right-sided-powerBingham}) exhibits a symmetric pattern with respect to girdle-type alternatives ($\tau<0$) and bipolar alternatives ($\tau>0$). This power, unlike the powers of the specified-$\thetab$ tests in the previous section, converges to~$\alpha$ as~$p$ diverges to infinity, which materializes the fact that, for larger dimensions, the Bingham test severely suffers (even asymptotically) from the unspecification of~$\thetab$.  Note also that since the Bingham test is invariant with respect to rotations, its limiting power naturally does not depend on the location parameter $\thetab$ under the alternative.

\begin{figure}
\centering
\includegraphics[width=\textwidth,height=70mm]{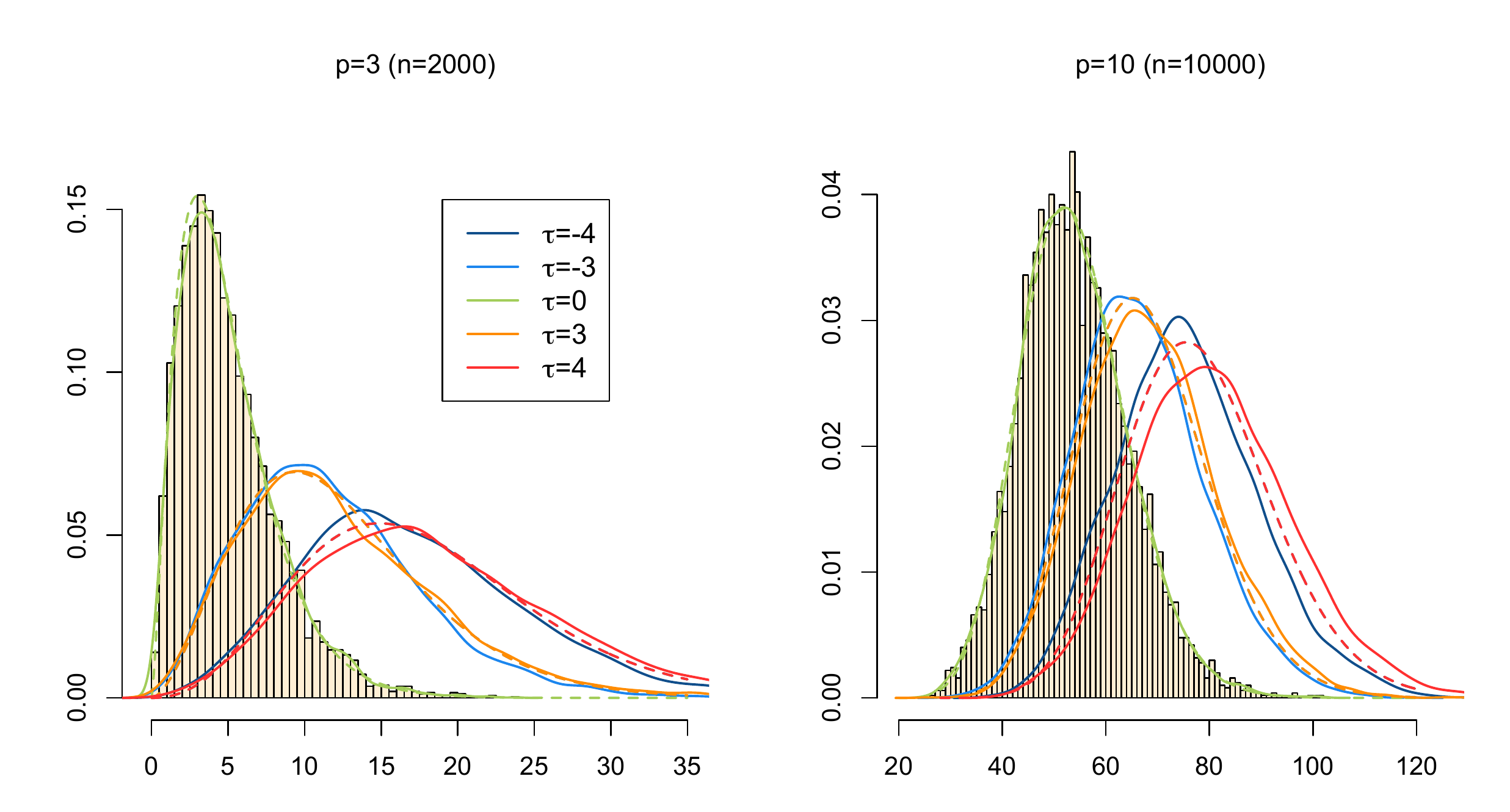}
\caption{
(Left:) Plots of the kernel density estimates (solid curves) of the values of the Bingham test statistic~$Q\n$ in~(\ref{Binghamtest}) obtained from~$M=5\:\!000$ independent random samples of size~$n=2\:\!000$ from the Watson distribution with location~$\thetab=(1,0,\ldots,0)'\in\R^p$ and concentration~$\kappa=\tau p/\sqrt{n}$, with~$\tau=-4,-3,0,3,4$ and with~$p=3$; for~$\tau=0$, histograms of the values of~$Q\n$ are shown. The density of the corresponding asymptotic distributions in~(\ref{Binghamtestncp}), which do not depend on the sign of~$\tau$, are also plotted (dashed curves). (Right:) The corresponding results for~$p=10$ and~$n=10\:\!000$. 
}  
\label{Fig2} 
\end{figure}

We conducted the following simulation exercise to check the validity of the asymptotic results of this section. In dimension~$p=3$, we generated $5\:\!000$ mutually independent random samples of size~$n=2\:\!000$ from the Watson distribution with location~$\thetab=(1,0,\ldots,0)'\in\R^{p}$ and concentration~$\kappa_n=\tau p/\sqrt{n}$, for $\tau=-4,-3,0,3,4$; see Section~\ref{sec:axial}. We did the same in dimension~$p=10$, with sample size~$n=10\:\!000$. For both dimensions~$p$, Figure~\ref{Fig2} reports kernel density estimates of the resulting values of the Bingham test statistic~$Q\n$. They perfectly match with the corresponding asymptotic distribution in~(\ref{Binghamtestncp}). The results also confirm the two-sided nature of the Bingham test, that, irrespective of~$\tau_0$, asymptotically behaves in the exact same way under~$\tau=\pm\tau_0$.

Our results indicate that the Bingham test shows non-trivial asymptotic powers under the contiguous alternatives identified in the previous section. 
\vspace{-.5mm}
However, a key point is the following: if~$(\taub_n)\to\taub$ in the 
\vspace{-.9mm}
LAN result of Theorem~\ref{TheorLANunspec}, then, under~${\rm P}_{\varthetab_n,f}\n$ with~$\varthetab_n=(p/\sqrt{n})^{1/2}\taub_n$,~$\Deltab\n$ is asymptotically normal with mean~${\bf s}_\taub=\Gamb {\rm vec}(\taub\taub')$ and covariance matrix~$\Gamb$, so that the sequence of asymptotic experiments at hand does converge to a Gaussian shift  experiment ($\Deltab\sim\mathcal{N}({\bf s}_\taub,\Gamb)$) involving a \emph{constrained} shift~${\bf s}_\taub$. As a result, the Bingham test is not Le Cam optimal
 for the considered problem: this test, which would rather be Le Cam
\vspace{-.6mm}
 optimal (more precisely, locally asymptotically maximin) for an unconstrained shift~${\bf s}\in\R^{p^2}$, is here ``wasting" power against multi-spiked alternatives that are incompatible with the present single-spiked axial model. This is in line with the fact that the Bingham test, which rejects the null hypothesis of uniformity when the sample variance of the eigenvalues~$\hat{\lambda}_{n1},\ldots,\hat{\lambda}_{np}$ of~$\Sb_n$ is too large, uses these eigenvalues in a permutation-invariant way 
\vspace{-.5mm}
(in the considered single-spiked models, it would be more natural to consider specifically~$\hat{\lambda}_{n1}$ and/or~$\hat{\lambda}_{np}$ to detect possible deviations from uniformity).


\section{The unspecified location case: single-spiked tests}
\label{sec:eigentest}

A natural question is then: how to construct a test that is more powerful than the Bingham test? We now describe two constructions that actually lead to the same test(s). Focusing again at first on the one-sided problem involving the bipolar alternatives, we saw in Section~\ref{sec:TestSpeci} that, in the specified location case, Le Cam optimal tests of uniformity reject~$\mathcal{H}_0\n:\{{\rm P}_0\n\}$ in favor of $\cup_{\kappa>0}\cup_{f\in\mathcal{F}} \{{\rm P}_{\thetab,\kappa,f}\n\}$ for large values of~$\Delta_{\thetab}\n=
\sqrt{n}
\big(
p
\,
\thetab'
\Sb_n
\thetab
-
1
\big)$. 
In the unspecified location case, it is
\vspace{-.7mm}
 then natural, following~\cite{Dav1977,Dav1987,Dav2002}, to consider the test,~$\phi_+^{(n)}$ say, rejecting the null hypothesis of uniformity at asymptotic level~$\alpha$ when
\begin{equation}
\label{Teststatlambda1}	
T_+^{(n)}
:=
\sup_{\thetab\in\mathcal{S}^{p-1}}
\Delta_{\thetab}\n
=
\sqrt{n}
\big(
p
\hat{\lambda}_{n1}
-
1
\big)
>
c_{p,\alpha,+}
,
\end{equation}
where~$\hat{\lambda}_{n1}$ still denotes the largest eigenvalue of~$\Sb_n$ and~$c_{p,\alpha,+}$ is such that this test has asymptotic size~$\alpha$ under the null hypothesis. 

A similar rationale yields 
\vspace{-.5mm}
natural tests for the other one-sided problem and for the two-sided problem: since 
\vspace{-.7mm}
Le Cam optimal tests of uniformity reject~$\mathcal{H}_0\n:\{{\rm P}_0\n\}$ in favor of $\cup_{\kappa<0}\cup_{f\in\mathcal{F}} \{{\rm P}_{\thetab,\kappa,f}\n\}$ for
\vspace{-.9mm}
 large values of~$-\Delta_{\thetab}\n=
\sqrt{n}
\big(
p
\,
\thetab'
\Sb_n
\thetab
-
1
\big)$, the resulting unspecified-$\thetab$ test,~$\phi_-^{(n)}$ say, will reject the null hypothesis of uniformity at asymptotic level~$\alpha$ when
\begin{equation}
\label{Teststatlambdap}	
T_-^{(n)}
:=
\sup_{\thetab\in\mathcal{S}^{p-1}}
(-\Delta_{\thetab}\n)
=
-
\sqrt{n}
\big(
p
\hat{\lambda}_{np}
-
1
\big)
>
c_{p,\alpha,-}
,
\end{equation}
where~$c_{p,\alpha,-}$ is such that this test has asymptotic size~$\alpha$ under the null hypothesis. Finally, since
\vspace{-.4mm}
 Le Cam optimal tests of uniformity reject~$\mathcal{H}_0\n:\{{\rm P}_0\n\}$ in favor of $\cup_{\kappa\neq 0}\cup_{f\in\mathcal{F}} \{{\rm P}_{\thetab,\kappa,f}\n\}$ for large values of~$|\Delta_{\thetab}\n|=
\sqrt{n}
\big|
p
\,
\thetab'
\Sb_n
\thetab
-
1
\big|$, the resulting unspecified-$\thetab$ test,~$\phi_\pm^{(n)}$ say, will reject the null hypothesis of uniformity at asymptotic level~$\alpha$ when
\begin{equation}
\label{Teststatlambda1p}	
T_{\pm}^{(n)}
:=
\sup_{\thetab\in\mathcal{S}^{p-1}}
|\Delta_{\thetab}\n|
=
\sqrt{n}
\max\big(
|p\hat{\lambda}_{n1}-1|
,
|p\hat{\lambda}_{np}-1|
\big)
>
c_{p,\alpha,\pm}
,
\end{equation}
where~$c_{p,\alpha,\pm}$ is still such that this test has asymptotic size~$\alpha$ under the null hypothesis of uniformity. 

Another rationale for considering the above tests is the following. For the sake of brevity, let us focus on the one-sided problem involving the bipolar alternatives, that is, the ones associated with~$\kappa>0$. A natural idea to obtain an unspecified-$\thetab$ test is to replace~$\thetab$ in the corresponding optimal specified-$\thetab$ test~$\phi_{\thetab+}^{(n)}$ with an estimator~$\hat{\thetab}_n$.
\vspace{-.8mm}
 Now, under~${\rm P}\n_{\thetab,\kappa,f}$, we have~${\rm E}[\Xb_{n1}]={\bf 0}$ and
$$
{\rm E}[\Xb_{n1}\Xb_{n1}']
=
g_f(\kappa) 
\thetab\thetab' 
+
\frac{1-g_f(\kappa)}{p-1} 
(\mathbf{I}_p-\thetab\thetab')
$$
(see, e.g., Lemma~B.3(i) in \citealp{PaiVer17a}), with
$$
g_f(\kappa)
:=
{\rm E}_{\thetab,\kappa,f}[(\Xb_{n1}'\thetab)^2]
=
c_{p,\kappa,f} \int_{-1}^1 (1-s^2)^{(p-3)/2} s^2 f(\kappa s^2)\,ds
.
$$
It is easy to check that, for any~$f\in\mathcal{F}$, the function~$\kappa\mapsto g_f(\kappa)$ is differentiable at~$0$, with derivative~$g_f'(0)={\rm Var}\n_{0}[(\Xb_{n1}'\thetab)^2]>0$, where~${\rm Var}\n_{0}$ still denotes variance under~${\rm P}\n_0$. Consequently, for~$\kappa>0$ small, we have 
$
g_f(\kappa)
> 
g_f(0)
=
1/p
$,
so that~$\thetab$ is, up to an unimportant sign (recall that only the pair~$\{\pm\thetab\}$ is identifiable), the leading unit eigenvector of~${\rm E}[\Xb_{n1}\Xb_{n1}']$ (for many functions~$f$, including the Watson one~$f(z)=\exp(z)$, this remains true for any~$\kappa>0$). 
Therefore, 
\vspace{-.3mm}
a moment estimator of~$\thetab$ is the leading eigenvector~$\hat{\thetab}_n$ of~$\Sb_n=\frac{1}{n} \sum_{i=1}^n \Xb_{ni}\Xb_{ni}'$. 
\vspace{-.7mm}
Note that in the Watson parametric submodel~$\{{\rm P}\n_{\thetab,\kappa,\exp}:\kappa>0\}$, this estimator~$\hat{\thetab}_n$ is also the MLE of~$\thetab$. The resulting test then rejects the null hypothesis of uniformity for large values of 
$$
\Delta\n_{\hat{\thetab}_n}
:=
\sqrt{n}
\big(
p
\hat{\thetab}_n'
\Sb_n
\hat{\thetab}_n
-
1
\big)
=
T_+^{(n)}
,
$$
hence coincides with the test~$\phi_+\n$ in~(\ref{Teststatlambda1}).  
\vspace{-.7mm}
A similar reasoning for the other one-sided problem 
leads to the test~$\phi_{-}^{(n)}$.

The critical values in~(\ref{Teststatlambda1})--(\ref{Teststatlambda1p}) above can of course be obtained from the asymptotic distribution of the corresponding test statistics under the null hypothesis. 
\vspace{-.4mm}
For~$p=3$, the asymptotic null distributions 
\vspace{-.9mm}
of~$T\n_+$ and~$T\n_-$ were obtained in \cite{AndSte1972}, where the corresponding 
\vspace{-.7mm}
one-sided tests~$\phi\n_+$ and~$\phi\n_-$ were first proposed. We extend their result to the two-sided test statistic~$T\n_\pm$ and, more importantly, to the non-null case. The key to do so is the following result.

%
%

\begin{Theor}
\label{thlambda1asympt}
Fix~$p\in\{2,3,\ldots\}$ and~$f\in\mathcal{F}$. Let~$\Zb$ be a $p\times p$ random matrix such that~${\rm vec}\,\Zb\sim\mathcal{N}
(
{\bf 0},
{\bf V}_p
)
$, with
$
{\bf V}_p
=
(p/(p+2))
(\mathbf{I}_{p^2}+{\bf K}_p) 
- 
(2/(p+2))
 \mathbf{J}_{p}
.
$
Then, 
(i) under~${\rm P}_0\n$, 
$$
\bigg(
\begin{array}{c}
\sqrt{n} (p\hat{\lambda}_{n1}-1) \\
\sqrt{n} (p\hat{\lambda}_{np}-1) 
\end{array}
\bigg)
\stackrel{\mathcal{D}}{\to}
\bigg(
\begin{array}{c}
L^{\rm max}_{p} \\
L^{\rm min}_{p} 
\end{array}
\bigg)
,
$$
where~$L_{p,{\rm max}}$ (resp.,~$L_{p,{\rm min}}$) is the largest (resp., smallest) eigenvalue of~${\Zb}$; 
\linebreak
(ii) under ${\rm P}\n_{\thetab,\kappa_n,f}$, where $\kappa_n=\tau_n p/\sqrt{n}$ involves a real sequence~$(\tau_n)$ converging to~$\tau$,
$$
\bigg(
\begin{array}{c}
\sqrt{n} (p\hat{\lambda}_{n1}-1) \\
\sqrt{n} (p\hat{\lambda}_{np}-1) 
\end{array}
\bigg)
\stackrel{\mathcal{D}}{\to}
\bigg(
\begin{array}{c}
L^{\rm max}_{p,\tau} \\
L^{\rm min}_{p,\tau} 
\end{array}
\bigg)
,
$$
where~$L^{\rm max}_{p,\tau}$ (resp.,~$L^{\rm min}_{p,\tau}$) is the largest (resp., smallest) eigenvalue of~${\Zb}_\tau:=\Zb+(2\tau/(p+2))\Wb_\tau$, with
${\Wb}_\tau:={\rm diag}(p-1,-1,\ldots,-1)$ for~$\tau\geq 0$ and~${\Wb}_\tau:={\rm diag}(-1,\ldots,-1,p-1)$ for~$\tau<0$. 
\end{Theor}



A direct consequence of Theorem~\ref{thlambda1asympt}(i) is that simulations can be used to obtain arbitrarily precise estimates of the asymptotic critical values needed to implement the tests~$\phi_+^{(n)}$, $\phi_-^{(n)}$ and~$\phi_\pm^{(n)}$. For instance, 
\vspace{-.7mm}
the test~$\phi_+^{(n)}$ will reject the null hypothesis of uniformity at asymptotic level~$\alpha$ whenever~$T_+^{(n)}
=
\sqrt{n}
\big(
p
\hat{\lambda}_{n1}
-
1
\big)
>
\hat{c}^{(m)}_{p,\alpha,+}$, where~$\hat{c}^{(m)}_{p,\alpha,+}$ denotes the upper $\alpha$-quantile of~$m$ independent realizations of~$L_{p}^{\rm max}$. Interestingly, the following corollary shows that simulations can actually be 
\vspace{-.5mm}
avoided in dimensions~$p=2$ and~$p=3$, 
\vspace{-.7mm}
as the asymptotic null distribution of~$T\n_+$, $T\n_-$ and~$T\n_\pm$ can be explicitly determined for these values of~$p$ (the result for $T\n_+$ and~$T\n_-$ in dimension~$p=3$ in~(\ref{cdfp3a}) below agrees with the one from \citealp{AndSte1972}).

\begin{Corol}
\label{Corollambda1asympt}
(i) Under the null hypothesis of uniformity over~$\mathcal{S}^{1}$, the test statistics~$T_+\n$, $T_-\n$, and~$T_\pm\n$ converge weakly to~$L^{\rm max}_{2}$, 
where~$L^{\rm max}_{2}$ has cumulative distribution function
\begin{equation}
	\label{cdfp2}
\ell\mapsto (1-\exp(-\ell^2)) \, \mathbb{I}[\ell>0]
;
\end{equation}
%
(ii) under the null hypothesis of uniformity over~$\mathcal{S}^{2}$, the test statistics~$T_+\n$ and $T_-\n$ converge weakly to~$L^{\rm max}_{3}$,
where~$L^{\rm max}_{3}$ has cumulative distribution function
\begin{equation}
	\label{cdfp3a}
\ell
\mapsto
\big\{\Phi(\sqrt{5}\ell)+\Phi({\textstyle{\frac{\sqrt{5}\ell}{2}}})+3\Phi''({\textstyle{\frac{\sqrt{5}\ell}{2}}})-1\big\} \, \mathbb{I}[\ell>0]
,
\end{equation}
whereas the test statistic~$T_\pm\n$ converges weakly to~$L_3:=\max(L^{\rm max}_{3},-L^{\rm min}_{3})$, 
where~$L_{3}$ has cumulative distribution function
\begin{equation}
	\label{cdfp3b}
\ell
\mapsto
\big\{
2\Phi({\textstyle{\frac{\sqrt{5}\ell}{2}}})
+
6\Phi''({\textstyle{\frac{\sqrt{5}\ell}{2}}})
-
2\sqrt{3} \Phi''({\textstyle{\frac{\sqrt{5}\ell}{\sqrt{3}}}})-1
\big\} 
\, \mathbb{I}[\ell>0]
\end{equation}
$($here, $\Phi''$ is the second derivative of the standard normal distribution function~$\Phi)$. 
\end{Corol}

Writing~$\lambda_\ell(\Ab)$ for the $\ell$th largest eigenvalue of the $p\times p$ matrix~$\Ab$ and denoting as~$\stackrel{\mathcal{D}}{=}$ equality in distribution, Theorem~\ref{thlambda1asympt} entails that, under the null hypothesis, 
$$
T_+\n
\stackrel{\mathcal{D}}{\to}
L_p^{\rm max}
=
\lambda_1(\Zb)
\stackrel{\mathcal{D}}{=}
\lambda_1(-\Zb)
=
-\lambda_p(\Zb)
=
-L_p^{\rm min}
\stackrel{\mathcal{D}}{\leftarrow}
T_-\n
.
$$ 
This shows that, for any dimension~$p$, the test statistics~$T_+\n$ and $T_-\n$ share the same weak limit under the null hypothesis, which is confirmed in dimensions~$p=2,3$ by Corollary~\ref{Corollambda1asympt}. Maybe surprisingly, this corollary further implies that, for~$p=2$, the two-sided test statistic~$T_\pm\n$ has the same asymptotic null distribution as~$T_+\n$ and $T_-\n$. 
\vspace{-.4mm}
This can be explained as follows: since~$\Sb_n$ has trace one almost surely, its eigenvalues~$\hat{\lambda}_{n\ell}$, $\ell=1,\ldots,p$ sum up to one almost surely (incidentally, this implies that the quantities~$\sqrt{n}(p\hat{\lambda}_{n\ell}-1)$, $\ell=1,\ldots,p$, do not admit a joint density, not even asymptotically so, which makes the proof of Theorem~\ref{thlambda1asympt} rather challenging). For~$p=2$, it follows that~$T_+\n=T_-\n=T_\pm\n$ almost surely, which of course entails that these three tests statistics share the same weak limit, not only under the null hypothesis but under \emph{any} sequence of hypotheses. In line with this, $L_{2,\tau}^{\rm max}=-L_{2,\tau}^{\rm min}=\max(L_{2,\tau}^{\rm max},-L_p^{\rm min})$ almost surely for any~$\tau\in\R$. 

\begin{figure}
\centering
\includegraphics[width=.95\textwidth]{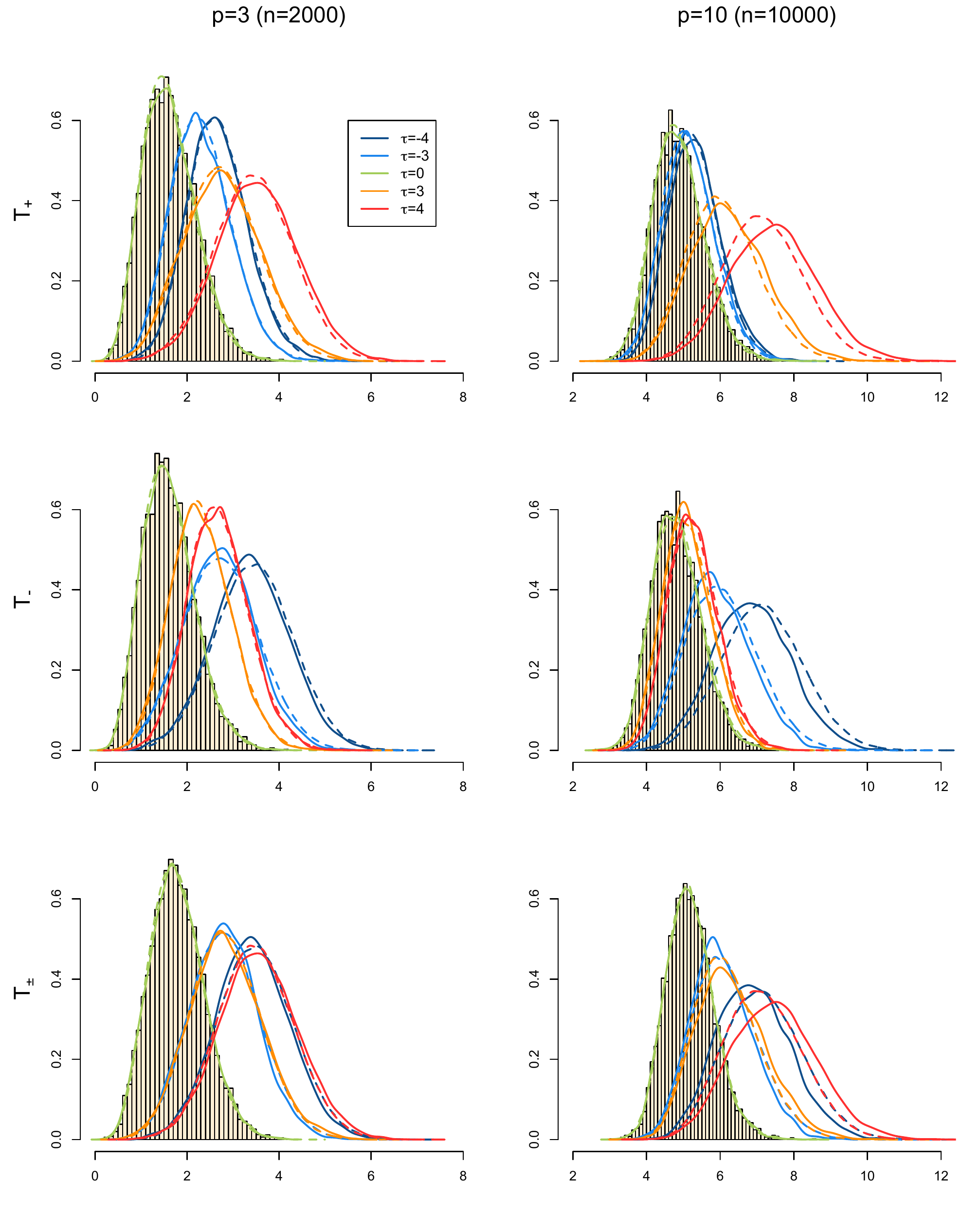}
\caption{
(Left:) Plots of 
\vspace{-.8mm}
the kernel density estimates (solid curves) of the values of~$T_+\n=\sqrt{n}(p\hat{\lambda}_{n1}-1)$ (top), $T_-\n=-\sqrt{n}(p\hat{\lambda}_{np}-1)$ (middle) or~$T_\pm\n=\max(T_+\n,T_-\n)$ (bottom), obtained from~$5\,000$ independent random samples of size~$n=2\,000$ from the Watson distribution with location~$\thetab=(1,0,\ldots,0)'\in\R^p$ and concentration~$\kappa=\tau p/\sqrt{n}$, with~$\tau=-4,-3,0,3,4$ and with~$p=3$; for~$\tau=0$, histograms of the corresponding test statistics are shown. The density of the corresponding asymptotic distributions  are also plotted (dashed curves). (Right:) The corresponding results for~$p=10$ and~$n=10\,000$. See Section~\ref{sec:eigentest} for details. 
}
\label{Fig3}
\end{figure}

To check the validity of Theorem~\ref{thlambda1asympt} and Corollary~\ref{Corollambda1asympt}, we conducted the following numerical exercises in dimensions~$p=3$ and~$p=10$. We generated $5\:\!000$ mutually independent random samples of size~$n=2\:\!000$ (for~$p=3$) and~$n=10\:\!000$ (for~$p=10$) from the Watson distribution with location~$\thetab=(1,0,\ldots,0)'\in\R^{p}$ and concentration~$\kappa_n=\tau p/\sqrt{n}$, for $\tau=-4,-3,0,3,4$. Figure~\ref{Fig3} plots kernel density estimates of the 
\vspace{-.4mm}
resulting values of~$T\n_+$, $T\n_-$ and~$T\n_\pm$, along with the densities of the corresponding asymptotic distributions; for~$p=3$ and~$\tau=0$, these densities are those associated with the distribution functions in~(\ref{cdfp3a})--(\ref{cdfp3b}), whereas, in all other cases, they are kernel density estimates obtained from~$10^6$ independent realizations of~$L^{\rm max}_{p,\tau}$, $-L^{\rm min}_{p,\tau}$, and $\max(L^{\rm max}_{p,\tau},-L^{\rm min}_{p,\tau})$, respectively; see Theorem~\ref{thlambda1asympt}. Clearly, 
\vspace{-.6mm}
the results support our asymptotic findings. It is seen that the one-sided test~$\phi_+\n$ not only shows power against the bipolar alternatives it is designed for (those associated with~$\tau>0$) but also against girdle-type ones (those associated with~$\tau<0$), which is actually desirable. 
\vspace{-.7mm}
The same can be said about the one-sided test~$\phi_-\n$, but 
\vspace{-.7mm}
each of these tests, of course, shows higher powers against the alternatives it was designed for. In contrast, the two-sided test~$\phi_\pm\n$ shows a symmetric power pattern for positive and negative values of~$\tau$.

\section{Finite-sample comparisons}
\label{sec:Simu}

In the previous sections, we conducted Monte-Carlo exercises in order to check correctness of our null and non-null asymptotic results, but it is of course of primary importance to compare the power behaviors of the various tests considered in this work.
\vspace{-.3mm}
 In this section, we 
\vspace{-1.3mm}
therefore study the finite-sample powers of the Bingham test~$\phi\n_{\rm Bing}$ and of the test~$\phi_+\n$ (we could similarly consider the tests~$\phi_+\n$ and~$\phi_\pm\n$), and %
\vspace{-.7mm}
we compare them with those of the optimal specified-$\thetab$ test~$\phi_{\thetab +}\n$. Our asymptotic results further allow us to complement these finite-sample comparisons with comparisons of the corresponding asymptotic powers.

We conducted the following Monte Carlo experiment. For any combination~$(n,p)$ of sample size~$n\in\{200,20\:\!000\}$ and dimension~$p\in\{3,10\}$, we generated collections of $2\:\!000$ independent random samples of size~$n$ from the Watson distribution on~$\mathcal{S}^{p-1}$ with location~$\thetab=(1,0,\ldots,0)'\in\R^p$ and concentration~$\kappa_n =\tau_\ell p/\sqrt{n}$, with~$\tau_\ell=0.8\ell$, $\ell=0,1,\ldots,5$. The value~$\ell=0$ corresponds to the null hypothesis of uniformity, whereas $\ell=1,\ldots,5$ provide
\vspace{-.4mm}
 increasingly severe bipolar alternatives. In each sample, we performed three tests at 
\vspace{-.9mm}
asymptotic level~$\alpha=5\%$, namely the specified-$\thetab$ test~$\phi_{\thetab+}^{(n)}$ in~(\ref{Testright}), the Bingham 
\vspace{-1.0mm}
test~$\phi\n_{\rm Bing}$ in~(\ref{Binghamtest}), and the test~$\phi\n_+$ in~(\ref{Teststatlambda1}); for~$p=3$, the asymptotic critical value for~$\phi_{+}^{(n)}$ was obtained from Corollary~\ref{Corollambda1asympt}(ii), whereas, for~$p=10$, an approximation of the corresponding critical value was obtained from~$10\:\!000$ independent realizations of~$L_p^{\rm max}$ in Theorem~\ref{thlambda1asympt}. 

\begin{figure}
\centering
\includegraphics[width=\textwidth]{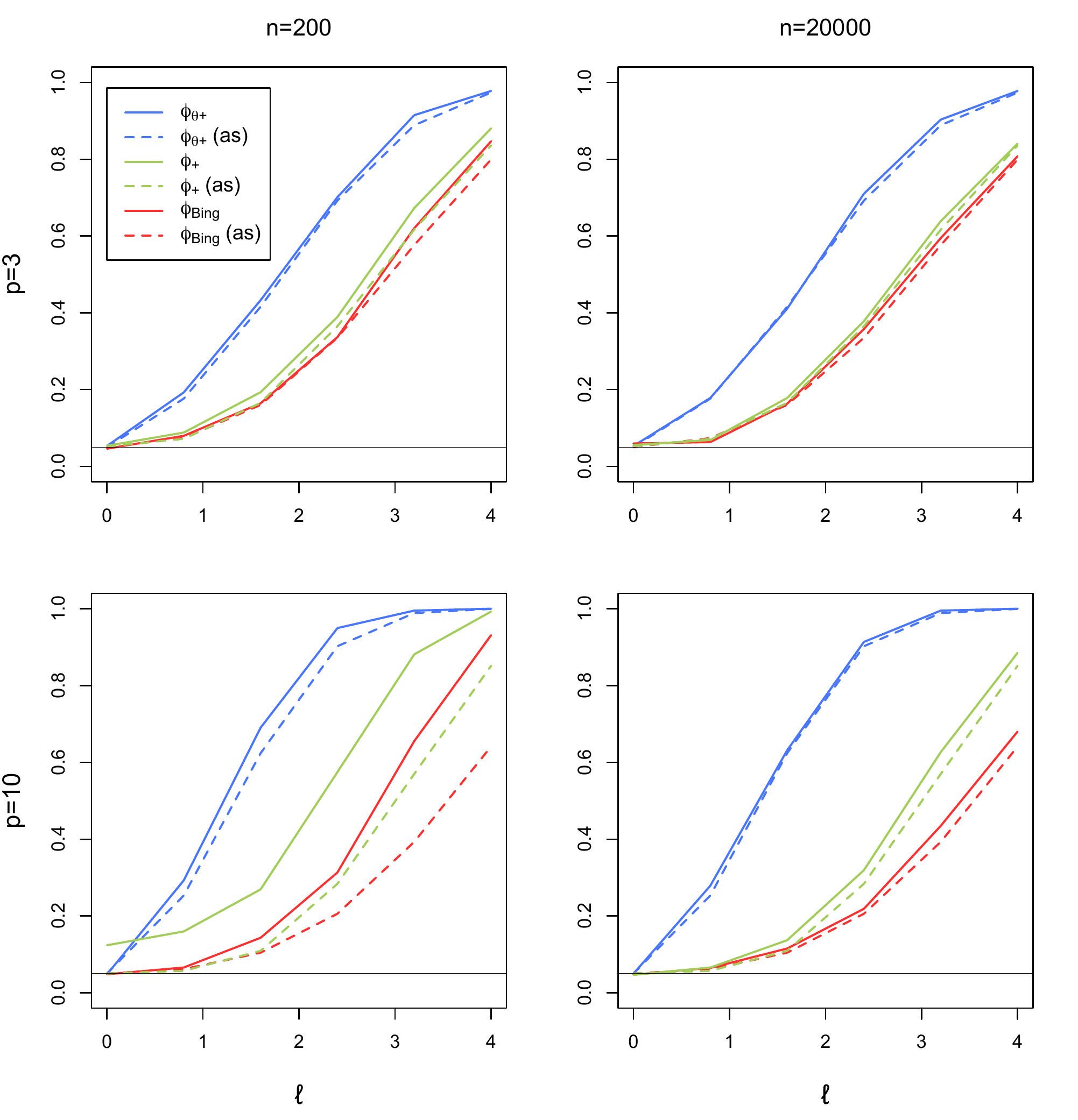}
\caption{
Rejection frequencies (solid curves) of three axial tests of uniformity over~$\mathcal{S}^{p-1}$ obtained from~$2\,000$ mutually independent random samples of size~$n$ from the Watson distribution with location~$\thetab=(1,0,\ldots,0)'\in\R^p$ and concentration~$\kappa=\tau_\ell p/\sqrt{n}$, with~$\ell=0$ (null hypothesis of uniformity) and~$\tau_\ell=0.8,1.6,2.4,3.2,4.0$ (increasingly
\vspace{-.7mm}
 severe bipolar alternatives). The tests
 %
 %
 considered are the test~$\phi_{\thetab+}^{(n)}$ in~(\ref{Testright}), the Bingham test~$\phi\n_{\rm Bing}$ 
 \vspace{-.3mm}
 in~(\ref{Binghamtest}), and the test~$\phi_+\n$ in~(\ref{Teststatlambda1}). The corresponding asymptotic powers (dashed curves) are also provided; see Section~\ref{sec:Simu} for details.
}
\label{Fig4}
\end{figure}

Figure~\ref{Fig4} 
\vspace{-.6mm}
shows the resulting empirical powers along with their theoretical asymptotic counterparts (for any given~$p$ and~$\tau_\ell$, the asymptotic power of~$\phi_{+}^{(n)}$ was obtained from~$10\,000$ independent copies of the random variable~$L_{p,\tau_\ell}^{\rm max}$ in Theorem~\ref{thlambda1asympt}). The results show that, 
\vspace{-.6mm}
as expected, the optimal specified-$\thetab$ test outperforms both unspecified-$\thetab$ tests. The test~$\phi\n_+$ dominates the Bingham test~$\phi\n_{\rm Bing}$ and this dominance, quite intuitively, increases with the dimension~$p$.  Clearly, rejection frequencies agree very well with our asymptotic results for large sample sizes.

%
%
%
%
%
%

%
%
\section{Final comments and perspectives for future research}
\label{sec:Final}

Practitioners often face axial data, which explains that statistical procedures for such data are presented in most directional statistics monographs and have been the topic of numerous research papers; we refer to \cite{Fish87}, \cite{MarJup2000}, \cite{LV17book}, and to the references therein. However, the non-null properties of axial tests of hypotheses are barely known, compared to those of their non-axial counterparts. Since this is in particular the case for axial tests of uniformity, we systematically studied in this work the asymptotic power of several axial tests of uniformity. In particular, we derived the asymptotic powers of the Bingham and Anderson tests under contiguous rotationally symmetric alternatives. Our results identify the underlying contiguity rate and allow for theoretical power comparisons. Throughout, our asymptotic findings were confirmed through simulations. Far from being of academic interest only, our results may be useful to practitioners who, when rejection occurs, will get some insight on the underlying distribution by combining the outcomes of the various tests of uniformity (that is, they will get hints on the single-spiked vs multi-spiked nature of the distribution, on its bipolar vs girdle nature, etc.)

The non-null asymptotic analysis conducted in this work essentially settles the low-dimensional case. Perspectives for future research therefore mainly relate to the high-dimensional framework. It is rather straightforward to extend the contiguity and LAN results in Theorems~\ref{TheorContig}--\ref{TheorLAN} to the case where the dimension~$p=p_n$ diverges to infinity with~$n$ at an arbitrary rate. However, in high dimensions, it is extremely challenging to derive the non-null asymptotic powers of the Bingham and Anderson tests under suitable local alternatives. For the Anderson tests, for instance, this is due to the fact that eigenvalues of sample covariance matrices suffer complicated phase transition phenomenons which, close to uniformity, results in a lack of consistency. These challenging questions are of course beyond the scope of the present work, hence are left for future research.


\appendix

\section{Proofs of Theorems~\ref{TheorContig}
to~\ref{TheorBinghamHD}}
\label{proofcontigth}

These proofs require the following preliminary result.

\begin{Lem}
\label{lemmeA1}
Let $g:\R\to\R$ be twice differentiable at 0, $(\kappa_n)$ be a sequence in~$\R_0$, and~$(p_n)$ be a sequence in~$\{2,3,\ldots\}$. Assume that~$\kappa_n$ is $o(p_n)$ as $n\to\infty$. Then,
\begin{eqnarray*}
R_n(g)
&\!\!:=\!\!&
c_{p_n}\int_{-1}^1(1-s^2)^{(p_n-3)/2}g(\kappa_n s^2)ds
\\[2mm]
&\!\!=\!\!&
g(0)+\frac{\kappa_n}{p_n}g'(0)+\frac{3\kappa_n^2}{2p_n(p_n+2)}g''(0)+o\left(\frac{\kappa_n^2}{p_n^2}\right)
,
\end{eqnarray*}
with~$c_p=1/\int_{-1}^1 (1-s^2)^{(p-3)/2}\,dt$.
\end{Lem}

{\sc Proof of lemma \ref{lemmeA1}.}
Recall that if~$\Xb$ is uniformly distributed over~$\mathcal{S}^{p_n-1}$, then, for any~$\thetab\in\mathcal{S}^{p_n-1}$, we have that~$\thetab'\Xb$ has density~$c_{p_n} (1-s^2)^{(p_n-3)/2}\mathbb{I}[|s|\leq 1]$, that is, the distribution of~$\thetab'\Xb$ is symmetric about zero and such that~$(\thetab'\Xb)^2\sim {\rm Beta}(1/2,(p_n-1)/2)$. Consequently, 
\begin{equation}
\label{foref1}
c_{p_n}
\int_{-1}^1 s^2 (1-s^2)^{(p_n-3)/2}  \,ds
=
{\rm E}[(\thetab'\Xb)^2]
=
\frac{1}{p_n}
,
\end{equation}
and
\begin{equation}
\label{foref2}
c_{p_n}
\int_{-1}^1 s^4 (1-s^2)^{(p_n-3)/2}  \,ds
=
{\rm E}[(\thetab'\Xb)^4]
=
\frac{3}{p_n(p_n+2)}
\cdot
\end{equation}
From~(\ref{foref1}), we can write
$$
R_n(g)-g(0)-\frac{\kappa_n}{p_n}g'(0)
=
c_{p_n}
\int_{-1}^1 (1-s^2)^{(p_n-3)/2} \left(g(\kappa_n s^2)-g(0)-\kappa_n s^2 g'(0)\right)  \,ds
.
$$
Letting~$t=|\kappa_n|^{1/2} s$ and using~(\ref{foref1})--(\ref{foref2}) then provides
$$
R_n(g)-g(0)-\frac{\kappa_n}{p_n}g'(0)
=
\frac{3\kappa_n^2}{p_n(p_n+2)}
\int_{-\infty}^{\infty} 
h_n(t) 
\left( 
\frac{g\big(\frac{\kappa_n}{|\kappa_n|}t^2\big)-g(0)-\frac{\kappa_n}{|\kappa_n|}t^2g'(0)}{t^4} 
\right)
\,dt
,
$$
or, equivalently,
\begin{eqnarray}
	\lefteqn{
\frac{R_n(g)-g(0)-\frac{\kappa_n}{p_n}g'(0)-\frac{3\kappa_n^2}{2p_n(p_n+2)}g''(0)}{\frac{3\kappa_n^2}{p_n(p_n+2)}}
}
\nonumber
\\[2mm]
& & 
\hspace{13mm} 
=
\int_{-\infty}^{\infty} 
h_n(t) 
\left( 
\frac{g\big(\frac{\kappa_n}{|\kappa_n|}t^2\big)-g(0)-\frac{\kappa_n}{|\kappa_n|}t^2g'(0)}{t^4} 
\right)
\,dt
-
\frac{1}{2}g''(0)
,
\label{prequeA1}
\end{eqnarray}
where~$h_n$ is defined through 
$$
t
\mapsto
h_n(t)
=
\frac{
\Big(\frac{t}{|\kappa_n|^{1/2}}\Big)^4
\Big(1-\Big(\frac{t}{|\kappa_n|^{1/2}}\Big)^2\Big)^{(p_n-3)/2}
\mathbb{I}\big[ |t|\leq |\kappa_n|^{1/2} \big]
}{
\int_{-|\kappa_n|^{1/2}}^{|\kappa_n|^{1/2}}
\Big(\frac{t}{|\kappa_n|^{1/2}}\Big)^4 
\Big(1-\Big(\frac{t}{|\kappa_n|^{1/2}}\Big)^2\Big)^{(p_n-3)/2}\,dt}
\cdot
$$
Since $\kappa_n$ is $o(p_n)$, it can be checked that the sequence~$(h_n)$ is an \emph{approximate \mbox{$\delta$-sequence}}, in the sense that
$
\int_{-\infty}^\infty h_n(t)\,dt=1 
$
for any~$n$ and
$
\int_{-\varepsilon}^\varepsilon h_n(t)\, dt \to 1
$
for any~$\varepsilon>0$. Hence, 
$$
\lim_{n\to \infty}
\int_{-\infty}^{\infty} 
h_n(t) 
\left( 
\frac{g\big(\frac{\kappa_n}{|\kappa_n|}t^2\big)-g(0)-\frac{\kappa_n}{|\kappa_n|}t^2g'(0)}{t^4} 
\right)
\,dt
=
\lim_{t\to 0} 
\frac{g\big(\frac{\kappa_n}{|\kappa_n|}t^2\big)-g(0)-\frac{\kappa_n}{|\kappa_n|}t^2g'(0)}{t^4}
,
$$
which, by using L'H\^{o}pital's rule, is equal to
$$
\lim_{t\to 0} \frac{ \frac{2\kappa_n}{|\kappa_n|}t g'\left(\frac{\kappa_n}{|\kappa_n|}t^2\right)-\frac{2\kappa_n}{|\kappa_n|}tg'(0)}{4t^3}
=
\frac{1}{2} 
\lim_{t\to 0} \frac{ g'\big(\frac{\kappa_n}{|\kappa_n|}t^2\big)-g'(0)}{\frac{\kappa_n}{|\kappa_n|} t^2} 
=
\frac{1}{2} g''(0)
.
$$
The result thus follows from~(\ref{prequeA1}). 
\cqfd
\vspace{3mm}


The proof of Theorems~\ref{TheorContig}--\ref{TheorLAN} actually only requires the particular case of lemma~\ref{lemmeA1} corresponding to~$p_n\equiv p$. We still presented this more general version of the lemma as it allows one to extend Theorems~\ref{TheorContig}--\ref{TheorLAN} to high-dimensional asymptotic scenarios where~$p_n$ would diverge to infinity with~$n$ (this would prove the claims in high dimensions provided in the last paragraph of Section~\ref{sec:axial}).
\vspace{2mm}
 
{\sc Proof of theorem \ref{TheorContig}.}
In this proof, all expectations and variances are taken under the null of uniformity~${\rm P}_{0}\n$ and all stochastic convergences and $o_{\rm P}$'s are as~$n\to\infty$ under~${\rm P}_{0}\n$. Throughout, 
we write $\ell_{f,k}:=(\log f)^k$ and~$E_{nk}:={\rm E}[ \ell_{f,k}(\kappa_n (\Xb_{ni}\pr\thetab)^2)]$. 

Consider the local log-likelihood ratio   
\begin{align*}
\Lambda_n
:=
\log \frac{d{\rm P}\n_{\thetab,\kappa_n,f}}{d{\rm P}\n_{0}}
&=
\sum_{i=1}^n 
\,
\log  \frac{c_{{p},\kappa_n,f} f(\kappa_n (\Xb_{ni}\pr \thetab)^2)}{c_{{p}}}\\
&= 
n 
\left(\log \frac{c_{{p},\kappa_n,f}}{c_{{p}}}+E_{n1}\right)
+
\sum_{i=1}^n 
\left(
\log f(\kappa_n (\Xb_{ni}\pr \thetab)^2) - E_{n1}
\right)\\
&=:
L_{n1}+L_{n2}.
\end{align*}
Lemma~\ref{lemmeA1} readily yields
\begin{align}
\log \frac{c_{{p},\kappa_n,f}}{c_{{p}}}
&=
- \log 
\left( 
c_{{p}}\int_{-1}^1 (1-s^2)^{({p}-3)/2} f(\kappa_n s^2)\,ds 
\right)
\nonumber\\[2mm]
&=
-\log 
\left(
1
+
\frac{\kappa_n}{{p}}
+
\frac{3\kappa_n^2}{2{p}({p}+2)}
f''(0) 
+
o(\kappa^2_n)
\right)
\nonumber\\[2mm]
&=
-\frac{\kappa_n}{{p}}
-
\frac{3\kappa_n^2}{2{p}({p}+2)}f''(0)
+
\frac{\kappa_n^2}{2{p}^2}
+
o(\kappa^2_n)
.
\label{expandlogconst}
\end{align}
Similarly,
\begin{align}
E_{n1}
&=
c_{{p}}
\int_{-1}^1 (1-s^2)^{({p}-3)/2}  \ell_{f,1}(\kappa_n s^2) \,ds 
\nonumber\\[2mm]
&=
\frac{\kappa_n}{{p}} \ell_{f,1}'(0)
+
\frac{3\kappa^2_n}{2{p}({p}+2)} \ell_{f,1}''(0)
+
o(\kappa^2_n)
\nonumber\\[2mm]
&=
\frac{\kappa_n}{{p}}+\frac{3\kappa^2_n}{2{p}({p}+2)} (f''(0)-1)
+
o(\kappa^2_n)
.
\label{expandEk}
\end{align}
Combining~(\ref{expandlogconst}) and~(\ref{expandEk}) provides
$$
L_{n1}
=
\frac{n\kappa_n^2}{2{p}^2}-\frac{3n\kappa^2_n}{2{p}({p}+2)} 
+
o(n\kappa^2_n)
=
-\frac{n\kappa_n^2({p}-1)}{{p}^2({p}+2)}
+
o(n\kappa^2_n)
.
$$

Turning to~$L_{n2}$, write
$$
L_{n2}
=
\sqrt{n V_n}
\,
\sum_{i=1}^n 
W_{ni}
:=
\sqrt{n V_n}
\,
\sum_{i=1}^n 
\frac{\log f(\kappa_n (\Xb_{ni}\pr \thetab)^2) - E_{n1}}{\sqrt{n V_n}}
,
$$
where we let
$
V_n
:=
{\rm Var}\big[ \log f(\kappa_n (\Xb_{ni}\pr\thetab)^2) \big]
.
$
First note that, since
\begin{align}
E_{n2}
&=
c_{{p}}
\int_{-1}^1 (1-s^2)^{({p}-3)/2}  \ell_{f,2}(\kappa_n s^2) \,ds 
\nonumber\\[2mm]
&=
\frac{\kappa_n}{{p}} \ell_{f,2}'(0)
+
\frac{3\kappa^2_n}{2{p}({p}+2)} \ell_{f,2}''(0)
+
o(\kappa^2_n)
\nonumber\\[2mm]
&=
\frac{3\kappa^2_n}{{p}({p}+2)}
+
o(\kappa^2_n)
,
\end{align}
we have 
\begin{equation}
\label{expandV}
nV_n
=
n\big(E_{n2} - E_{n1}^2\big)
=
\frac{2n\kappa_n^2({p}-1)}{{p}^2({p}+2)}
+
o(n\kappa^2_n)
,
\end{equation}
which leads to
\begin{equation}
\label{quasicont}
\Lambda_n
=
-\frac{n\kappa_n^2({p}-1)}{{p}^2({p}+2)}
+
\sqrt{\frac{2n\kappa_n^2({p}-1)}{{p}^2({p}+2)}
+
o(n\kappa^2_n)}
\,
\sum_{i=1}^n 
W_{ni}
+
o(n\kappa^2_n)
.
\end{equation}
Since $W_{ni}$, $i=1,\ldots,n$, are mutually independent with mean zero and variance~$1/n$, we obtain that
\begin{equation}
\label{L2norm}
{\rm E}
\big[
\Lambda_n^2
\big]
=
\left(
{\rm E}
[
\Lambda_n
]
\right)^2
+
{\rm Var}
[
\Lambda_n
]
=
\frac{n^2\kappa_n^4({p}-1)^2}{{p}^4({p}+2)^2}
+
o(n^2\kappa^4_n)
+
\frac{2n\kappa_n^2({p}-1)}{{p}^2({p}+2)}
+
o(n\kappa^2_n)
.
\end{equation}

If $\kappa_n=o(1/\sqrt{n})$, then~(\ref{L2norm}) implies that $\exp(\Lambda_n)=1+o_{\rm P}(1)$, so that Le Cam's first lemma yields that~${\rm P}\n_{\thetab,\kappa_n,f}$ and ${\rm P}\n_{0}$ are mutually contiguous. 

We may therefore assume that $\kappa_n=\tau_n/\sqrt{n}$, 
where~($\tau_n$) is~$O(1)$ but not~$o(1)$, or, equivalently, that $\kappa_n=\tau_n {p}/\sqrt{n}$ with the same sequence~($\tau_n$). Then, (\ref{quasicont}) rewrites
$$
\Lambda_n
=
-
\frac{({p}-1)\tau_n^2}{{p}+2}
+
\sqrt{
\frac{2({p}-1)\tau_n^2}{{p}+2} 
+
o(1)
}
\,
\sum_{i=1}^n
W_{ni}
+
o(1).
$$

Applying the Cauchy-Schwarz inequality and the Chebychev inequality, then using Lemma~\ref{lemmeA1} and~(\ref{expandV}), yields that, for some positive constant~$C$ and any~$\varepsilon>0$,
\begin{eqnarray*}
\lefteqn{
\sum_{i=1}^n 
{\rm E}[
W_{ni}^2 
\mathbb{I}[ |W_{ni}| >\varepsilon]]
}
\\[1mm]
& &
\hspace{1mm} 
\leq 
n
\sqrt{{\rm E}[ W_{ni}^4] {\rm P}[|W_{ni}|>\varepsilon]}
\leq
\frac{n}{\varepsilon}
\sqrt{{\rm E}[ W_{ni}^4] {\rm Var}[W_{ni}]}
=
\frac{1}{\varepsilon}
\sqrt{n{\rm E}[ W_{ni}^4]}
\leq
\frac{C\sqrt{n E_{n4}}}{\varepsilon nV_n}
\\[2mm]
& &
\hspace{1mm} 
=
\frac{C
\Big(
\frac{n\kappa_n}{{p}} \ell_{f,4}'(0)
+
\frac{3n\kappa_n^2}{2{p}({p}+2)} \ell_{f,4}''(0)
+
o(n\kappa^2_n)
\Big)^{1/2}
}{
\varepsilon
\Big(
\frac{2n\kappa_n^2({p}-1)}{{p}^2({p}+2)}
+
o(n\kappa^2_n)
\Big)
}
=
\frac{
o(\tau_n)
}{
\varepsilon
\Big(
\frac{2({p}-1)\tau_n^2}{{p}+2} 
+
o(\tau_n^2)
\Big)
}
=
o(1)
,
\end{eqnarray*}
where we have used the fact that $\ell_{f,4}'(0)=\ell_{f,4}''(0)=0$. This shows that $\sum_{i=1}^n W_{ni}$ satisfies the classical Levy-Lindeberg condition,  hence is asymptotically standard normal (as already mentioned, $W_{ni}$, $i=1,\ldots,n$, are mutually independent with mean zero and variance~$1/n$). For any subsequence~$(\exp(\Lambda_{n_m}))$ converging in distribution, the weak limit must then be
$
\exp(Y)
$,
with
$
Y
\sim
\mathcal{N}( 
- \eta
,
 2\eta 
)
,
$
where we let~$\eta:=((p-1)/(p+2)) \lim_{m\to\infty} \tau_{n_m}^2$. Mutual contiguity of ${\rm P}\n_{\thetab,\kappa_n,f}$ and ${\rm P}\n_{0}$ then follows from the fact that ${\rm P}[\exp(Y)=0]=0$ and ${\rm E}[\exp(Y)]=1$.
\cqfd
\vspace{3mm}



{\sc Proof of Theorem \ref{TheorLAN}.}
As in the proof of Theorem~\ref{TheorContig}, all expectations and variances in this proof are taken under the null of uniformity~${\rm P}_{0}\n$ and all stochastic convergences and $o_{\rm P}$'s are as~$n\to\infty$ under~${\rm P}_{0}\n$. 
%
%
Recall that we have obtained in the proof of Theorem~\ref{TheorContig} that 
\begin{eqnarray*}
\Lambda_n
&\!=\!&
-
\frac{({p}-1)\tau_n^2}{{p}+2}
+
\sqrt{
\frac{2({p}-1)\tau_n^2}{{p}+2} 
+
o(1)
}
\,
\sum_{i=1}^n
W_{ni}
+
o(1)
\\[2mm]
&\!=\!&
-
\frac{\tau_n^2}{2} \Gamma_p
+
|\tau_n| \Gamma_p^{1/2}
\,
\sum_{i=1}^n
W_{ni}
+
o_{\rm P}(1)
,
\end{eqnarray*}
where 
$
\sum_{i=1}^n W_{ni}
=
(1/\sqrt{nV_n})
\sum_{i=1}^n
(\log f(\kappa_n (\Xb_{ni}\pr \thetab)^2) - E_{n1})
$
is asymptotically standard normal. Since~$(\tau_n\Gamma_p^{1/2})$ is~$O(1)$, it is therefore sufficient to show that
\begin{equation}
d_n
:=
{\rm E}\bigg[
	 \bigg( \Gamma_p^{-1/2}\Delta_{\thetab}\n- \frac{\tau_n}{|\tau_n|}\sum_{i=1}^n W_{ni} \bigg)^2
	\bigg]
	=
	o(1)
.
\label{EM2a}
\end{equation}
To do so, write
\begin{eqnarray}
\lefteqn{	
\hspace{-0mm} 
\Gamma_p^{-1/2}\Delta_{\thetab}\n- \frac{\tau_n}{|\tau_n|}\sum_{i=1}^n W_{ni}
}
\nonumber
\\[2mm]
& &
\hspace{2mm} 
=
\frac{1}{\sqrt{nV_n}}
\sum_{i=1}^n
\bigg(
{p}\sqrt{V_n}
\Gamma_p^{-1/2} \Big((\Xb_{ni}\pr \thetab)^2-\frac{1}{{p}}\Big)
-
\frac{\tau_n}{|\tau_n|}
(\log f(\kappa_n (\Xb_{ni}\pr \thetab)^2) - E_{n1})
\bigg)
\nonumber
\\[2mm]
& &
\hspace{2mm} 
=:
\frac{M_n}{\sqrt{nV_n}} 
\cdot
\label{EM2b}
\end{eqnarray}

Then using the fact that
${\rm E}[(\Xb_{n1}\pr \thetab)^2]=1/{p}$
and
${\rm E}[(\Xb_{n1}\pr \thetab)^4]=3/({p}({p}+2))$ (see the proof of Lemma~\ref{lemmeA1}) provides~${\rm Var}[(\Xb_{n1}\pr \thetab)^2]=\Gamma_p/{p}^2$, we obtain
\begin{align*}
{\rm E}\big[M_n^2\big]&=
n
{\rm E}
\left[
\left(
{p}\sqrt{V_n}
\Gamma_p^{-1/2}
\left((\Xb_{n1}\pr \thetab)^2-\frac{1}{{p}}\right)
-
\frac{\tau_n}{|\tau_n|}
(\log f(\kappa_n (\Xb_{n1}\pr \thetab)^2) - E_{n1})
\right)^2
\right]
\\[2mm]
&=
2nV_n 
-
2n{p} \sqrt{V_n} \Gamma_p^{-1/2}
\frac{\tau_n}{|\tau_n|}
{\rm E}
\left[
\left((\Xb_{n1}\pr \thetab)^2-\frac{1}{{p}}\right)
(\log f(\kappa_n (\Xb_{n1}\pr \thetab)^2) - E_{n1})
\right]
,
\end{align*}
which, letting $g(x):=x\log f(x)$, rewrites
\begin{eqnarray}
{\rm E}\big[M_n^2\big]
&\!\!=\!\!&
2nV_n
-
2n{p} \sqrt{V_n} \Gamma_p^{-1/2}
\frac{\tau_n}{|\tau_n|}
\Big(
{\rm E}
\big[
(\Xb_{ni}\pr \thetab)^2 \log f(\kappa_n (\Xb_{n1}\pr \thetab)^2) 
\big]
- \frac{E_{n1}}{{p}}
\Big)
\nonumber
\\[2mm]
&\!\!=\!\!&
2nV_n
-
2 \sqrt{nV_n} \Gamma_p^{-1/2}
\frac{\tau_n}{|\tau_n|}
\Big(
\frac{\sqrt{n}{p}}{\kappa_n} 
{\rm E}
\big[
g(\kappa_n (\Xb_{n1}\pr \thetab)^2) 
\big]
-
\sqrt{n} E_{n1}
\Big)
.
\label{EM2}
\end{eqnarray}
Lemma~\ref{lemmeA1} provides
$$
{\rm E}[
g(\kappa_n (\Xb_{n1}\pr \thetab)^2) 
]
=
c_{{p}}
\int_{-1}^1 (1-s^2)^{({p}-3)/2}  g(\kappa_n s^2) \,ds 
=
\frac{3\kappa^2_n}{{p}({p}+2)} 
+
o(\kappa^2_n)
.
$$
Using this jointly with~(\ref{expandEk}) and~(\ref{expandV}), it follows from~(\ref{EM2b})--(\ref{EM2}) that 
\begin{eqnarray*}
d_n
&\!\!=\!\!&
2
-
\frac{
2\tau_n
\Big(
\frac{3\sqrt{n}\kappa_n}{{p}+2} 
+
o\left(\frac{\sqrt{n}\kappa_n}{{p}} \right)
-
\Big\{
\frac{\sqrt{n}\kappa_n}{{p}}+\frac{3\sqrt{n}\kappa^2_n}{2{p}({p}+2)} (f''(0)-1)
+
o(\sqrt{n}\kappa^2_n)
\Big\}
\Big)
}{
\sqrt{\Gamma_p}|\tau_n|
\left(
\frac{n\kappa_n^2\Gamma_p}{{p}^2}
+
o(n\kappa^2_n)
\right)^{1/2}
}
\\[2mm]
&\!\!=\!\!&
2
-
\frac{
2\tau_n
\big(
\Gamma_p \tau_n
+
o(1)
\big)
}{
\sqrt{\Gamma_p}|\tau_n|
\big(
\Gamma_p\tau_n^2
+
o(1)
\big)^{1/2}
}
=
o(1)
,
\end{eqnarray*}
as was to be shown. 
\cqfd
\vspace{3mm}
 

{\sc Proof of Theorem \ref{TheorLANunspec}.}
First note that by proceeding strictly along the same lines as in the proof of Theorem~\ref{TheorLAN}, one can show that, for any sequence~$(\thetab_n)$ in~$\mathcal{S}^{p-1}$, and any real sequence~$(\tau_n)$ that is~$O(1)$ but not~$o(1)$,  
\begin{equation}
\label{LANthetan}
\log \frac{d{\rm P}\n_{\thetab_n,\kappa_n,f}}{d{\rm P}\n_{0}} 
=
\tau_n
\Delta_{\thetab_n}\n
-
\frac{\tau_n^2}{2} \Gamma_p
+
o_{\rm P}(1)
,
\end{equation}
under~${\rm P}_0\n$, with~$\kappa_n=\tau_n p/\sqrt{n}$ (that is, the stochastic second-order expansion in~(\ref{LAN}) remains true if one replaces the fixed location~$\thetab$ with an~$n$-dependent value~$\thetab_n$). Therefore, noting that the parameter value~$\varthetab_n=(p/\sqrt{n})^{1/2}\taub_n$ in the statement of the theorem corresponds to~$\thetab_n=\taub_n/\|\taub_n\|$ and~$\kappa_n=p\|\taub_n\|^2/\sqrt{n}$, we have
$$
\log \frac{d{\rm P}\n_{\varthetab_n,f}}{d{\rm P}\n_{0}} 
=
\log \frac{d{\rm P}\n_{\taub_n/\|\taub_n\|,p\|\taub_n\|^2/\sqrt{n},f}}{d{\rm P}\n_{0}}
=
\|\taub_n\|^2
\Delta_{\taub_n/\|\taub_n\|}\n
-
\frac{\|\taub_n\|^4}{2} \Gamma_p
+
o_{\rm P}(1)
$$  
under~${\rm P}_0\n$. Now, since the central sequence in Theorem~\ref{TheorLAN} rewrites
\begin{eqnarray}
\lefteqn{
\Delta_{\thetab}\n
=
\frac{p}{\sqrt{n}}\sum_{i=1}^n \Big\{(\Xb_{ni}\pr \thetab)^2-\frac{1}{p}\Big\}
=
p\sqrt{n}
\,
\thetab'
\bigg(
\Sb_n
-\frac{1}{p}\mathbf{I}_{p}
\bigg)
\thetab
}
\nonumber
\\[2mm]
& & 
\hspace{18mm} 
=
p\sqrt{n}
\,
(\thetab\otimes \thetab)'
{\rm vec}\bigg(
\Sb_n-\frac{1}{p} \mathbf{I}_p
\bigg)
=
({\rm vec}(\thetab\thetab'))'
\Deltab\n
,
\label{centralrewr}
\end{eqnarray}
this rewrites
$$
\log \frac{d{\rm P}\n_{\varthetab_n,f}}{d{\rm P}\n_{0}} 
=
({\rm vec}(\taub_n\taub_n'))'
\Deltab\n
-
\frac{\|\taub_n\|^4}{2} \Gamma_p
+
o_{\rm P}(1)
$$
under~${\rm P}_0\n$. By using the identities~$({\rm vec}\,\Ab)'({\rm vec}\,\Bb)={\rm tr}[\Ab'\Bb]$ and ${\bf K}_p({\rm vec}\,\Ab)={\rm vec}(\Ab')$, straightforward calculations show that~$\|\taub_n\|^4 \Gamma_p=({\rm vec}(\taub_n\taub_n'))'\Gamb {\rm vec}(\taub_n\taub_n')$, which establishes~(\ref{LANvec}). Since the asymptotic normality result readily follows from the multivariate central limit theorem, the theorem is proved. 
\cqfd
\vspace{3mm}

%
%
%
%
%
%
%
%
%

{\sc Proof of Theorem~\ref{TheorBinghamHD}.}
Denoting as ${\rm E}_0\n$ and~${\rm Var}_0\n$ expectation and variance under~${\rm P}_0\n$, one has (see~(\ref{centralrewr}))
\begin{eqnarray*}
\lim_{n\to\infty}
{\rm E}_0\n
\big[
\Delta_{\thetab}\n
\Deltab\n
\big]
&\!=\!&
\lim_{n\to\infty}
{\rm E}_0\n
\Bigg[
\Deltab\n
\bigg(
p\sqrt{n}
\,
(\thetab\otimes \thetab)'
{\rm vec}\bigg(
\Sb_n
-\frac{1}{p}\mathbf{I}_{p}
\bigg)
\bigg)\pr 
\Bigg]
\\[2mm]
&\!=\!&
\lim_{n\to\infty}
{\rm Var}_0\n
\big[
\Deltab\n
\big]
(\thetab\otimes \thetab)
=
\Gamb
(\thetab\otimes \thetab)
.
\end{eqnarray*}
Therefore, Le Cam's third lemma implies that, under ${\rm P}\n_{\thetab,\kappa_n,f}$ with $\kappa_n=\tau p/\sqrt{n}$, ($\tau\neq 0$), $\Deltab\n$ is asymptotically normal with mean vector~$\tau\Gamb
(\thetab\otimes \thetab)$ and covariance matrix~$\Gamb$, so that, under the same sequence of hypotheses, 
$$
Q\n
=
(\Deltab\n)' {\pmb\Gamma}_p^- \Deltab\n
\stackrel{\mathcal{D}}{\to}
\chi^2_{d_p}(\delta)
,
$$
with
$$
\delta
=
\tau^2 
(\thetab\otimes \thetab)'
\Gamb{\pmb\Gamma}_p^-\Gamb
(\thetab\otimes \thetab)
=
\frac{p\tau^2}{p+2}\Big( 2-\frac{2}{p}\Big)
=
\frac{2(p-1)\tau^2}{p+2}
\cdot
$$
The asymptotic power in~(\ref{eq:right-sided-powerBingham}) readily follows. 
\cqfd
\vspace{3mm}


\section{Proofs of Theorem~\ref{thlambda1asympt} and Corollary~\ref{Corollambda1asympt}}
\label{proofeigentests}

{\sc Proof of Theorem~\ref{thlambda1asympt}.}
Since ${\rm E}[\Xb_{n1}\Xb_{n1}']=(1/p)\mathbf{I}_p$ and ${\rm E}[{\rm vec}(\Xb_{n1}\Xb_{n1}') 
({\rm vec}(\Xb_{n1}\Xb_{n1}'))']
\linebreak
=
1/(p(p+2))(\mathbf{I}_{p^2}+{\bf K}_p+{\bf J}_p)$ under~${\rm P}_0\n$ (see, e.g., Lemma~A.2 in \citealp{PaiVer2016}), the multivariate central limit theorem yields
$$
\sqrt{n} 
\, 
{\rm vec}\bigg(
\Sb_n-\frac{1}{p} \mathbf{I}_p
\bigg)
\stackrel{\mathcal{D}}{\to}
\mathcal{N}
\bigg(
{\bf 0}
,
\frac{1}{p^2}{\bf V}_p
\bigg)
.
$$
Now, by using~(\ref{centralrewr}), we obtain that, under~${\rm P}_0\n$, 
\begin{eqnarray*}
	\lefteqn{
\hspace{-13mm} 
{\rm E}\bigg[
\sqrt{n} 
\, 
{\rm vec}\bigg(
\Sb_n-\frac{1}{p} \mathbf{I}_p
\bigg)
\Delta_{\thetab}\n
\bigg]
=
np
{\rm E}\bigg[
{\rm vec}\bigg(
\Sb_n-\frac{1}{p} \mathbf{I}_p
\bigg)
{\rm vec}'\bigg(
\Sb_n-\frac{1}{p} \mathbf{I}_p
\bigg)
\bigg]
{\rm vec}(\thetab\thetab')
}
\\[2mm]
& & 
\hspace{3mm} 
=
p
\bigg(
\frac{1}{p^2} {\bf V}_p
\bigg)
{\rm vec}(\thetab\thetab')
=
\frac{2}{p+2} {\rm vec}(\thetab\thetab')
 - \frac{2}{p(p+2)} {\rm vec}(\mathbf{I}_p)
,
\end{eqnarray*}
so that Le Cam's third lemma shows that, under ${\rm P}\n_{\thetab,\kappa_n,f}$, where $\kappa_n=\tau_n p/\sqrt{n}$ is based on a sequence~$(\tau_n)$ converging to~$\tau$, 
$$
\sqrt{n} 
\, 
{\rm vec}\bigg(
\Sb_n-\frac{1}{p} \mathbf{I}_p
\bigg)
\stackrel{\mathcal{D}}{\to}
\mathcal{N}
\bigg(
\frac{2\tau}{p+2} {\rm vec}(\thetab\thetab')
 - \frac{2\tau}{p(p+2)} {\rm vec}(\mathbf{I}_p)
,
\frac{1}{p^2}
{\bf V}_p
\bigg)
,
$$
which rewrites
\begin{equation}
	\label{TCLvic}
\sqrt{n} 
\, 
{\rm vec}\big(
\Sb_n
-
\Sigb_n
 \big)
\stackrel{\mathcal{D}}{\to}
\mathcal{N}
\bigg(
{\bf 0},
\frac{1}{p^2}
{\bf V}_p
\bigg)
,
\end{equation}
where
\begin{eqnarray*}
\Sigb_n
&\!\!\!:=\!\!\!&
\bigg(
\frac{1}{p} - \frac{2\tau}{\sqrt{n}p(p+2)} 
\bigg)
\mathbf{I}_p
+
\frac{2\tau}{\sqrt{n}(p+2)} \thetab\thetab'
\\[2mm]
&\!\!\!=\!\!\!&
\bigg(
\frac{1}{p} + \frac{2(p-1)\tau}{\sqrt{n}p(p+2)}
\bigg)
\thetab\thetab'
+
\bigg(
\frac{1}{p} - \frac{2\tau}{\sqrt{n}p(p+2)} 
\bigg)
(\mathbf{I}_p-\thetab\thetab')
.
\end{eqnarray*}
We need to consider the cases (a)~$\tau\geq 0$ and (b)~$\tau<0$ separately.  

In Case~(a), $\Sigb_n$ has eigenvalues 
\begin{equation}
	\label{eigenvic}
\lambda_{n1}
=
\frac{1}{p} + \frac{2(p-1)\tau}{\sqrt{n}p(p+2)}
\quad
\textrm{ and }
\quad
\lambda_{n2}
=
\ldots
=
\lambda_{np}
=
\frac{1}{p} - \frac{2\tau}{\sqrt{n}p(p+2)}
\cdot
\end{equation}
Fix then arbitrarily~$\thetab_2, \ldots, \thetab_p$ such that the $p \times p$ matrix~${\bf G}_p:=(\thetab,\thetab_2, \ldots, \thetab_p)$ is orthogonal. Letting~$\Lamb_n:={\rm diag}(\lambda_{n1},\lambda_{n2}, \ldots, \lambda_{np})$, 
we have $\Sigb_n {\bf G}_p={\bf G}_p\Lamb_n$, so that ${\bf G}_p$ is an eigenvectors matrix for $\Sigb_n$. Clearly, $\xi_{n1}:=\sqrt{n}p (\hat{\lambda}_{n1}-\lambda_{n1})$ is the largest root of the polynomial~$P_{n1}(h):={\rm det}({\sqrt n}p({\bf S}_n-\lambda_{n1} {\bf I}_p)- h {\bf I}_p)$, whereas~$\xi_{np}:=\sqrt{n}p (\hat{\lambda}_{np}-\lambda_{np})$ is the smallest root of~$P_{np}(h):={\rm det}({\sqrt n}p({\bf S}_n-\lambda_{np} {\bf I}_p)- h {\bf I}_p)$. Letting
\begin{equation}
	\label{defZn}
{\bf Z}_n
:= 
\sqrt{n}p({\bf G}_p' {\bf S}_n {\bf G}_p- \Lamb_n)
=
 {\bf G}_p' \sqrt{n}p({\bf S}_n- \Sigb_n) {\bf G}_p
 ,
\end{equation}
rewrite these polynomials as
$
P_{nj}(h)
=
{\rm det}(\sqrt{n}p{\bf G}_p' {\bf S}_n {\bf G}_p- \sqrt{n}p\lambda_{nj} {\bf I}_p - h {\bf I}_p) 
=
{\rm det}({\Zb}_n+\sqrt{n}p(\Lamb_n-\lambda_{nj} {\bf I}_p)- h {\bf I}_p)
,
$
$j=1,p$, which gives
$$
P_{n1}(h)
=
{\rm det}({\Zb}_n+{\rm diag}(0, -v_\tau, \ldots, -v_\tau)- h {\bf I}_p)
\nonumber  
$$
and
$$
P_{np}(h)
=
{\rm det}({\Zb}_n+{\rm diag}(v_\tau,0, \ldots,0)- h {\bf I}_p)
\nonumber  
,
$$
where we wrote
$$
v_\tau
:=
\sqrt{n}p
\bigg[
\frac{2(p-1)\tau}{\sqrt{n}p(p+2)}
+
\frac{2\tau}{\sqrt{n}p(p+2)}
\bigg]
=
\frac{2p\tau}{p+2}
\cdot
$$

Note that~(\ref{TCLvic}) readily implies that~${\rm vec}\,\Zb_n
=
\sqrt{n}p 
({\bf G}_p\otimes {\bf G}_p)' 
\, 
{\rm vec}\big(
\Sb_n
-
\Sigb_n
 \big)
$ converges weakly to~${\rm vec}\,\Zb\sim\mathcal{N}
(
{\bf 0},
{\bf V}_p
)
$. It readily follows that~$(\xi_{n1},\xi_{np})'$ converges weakly to~$(\xi_{1},\xi_{p})'$, where~$\xi_{1}$ is the largest root of the polynomial~${\rm det}({\Zb}+ {\rm diag}(0, -v_\tau, \ldots, -v_\tau)- h {\bf I}_p)$ and~$\xi_{p}$ is the smallest root of the polynomial~${\rm det}({\Zb}+ {\rm diag}(v_\tau,0, \ldots, 0)- h {\bf I}_p)$. This implies that 
\begin{eqnarray*}
\bigg(
\begin{array}{c}
\sqrt{n} (p\hat{\lambda}_{n1}-1) \\
\sqrt{n} (p\hat{\lambda}_{np}-1) 
\end{array}
\bigg)
&\! = \!&
\bigg(
\begin{array}{c}
\xi_{n1}\\
\xi_{np}
\end{array}
\bigg)
+
\frac{2\tau}{p+2}
\bigg(
\begin{array}{c}
p-1 \\
-1
\end{array}
\bigg)
\\[2mm]
&\! \stackrel{\mathcal{D}}{\to} \!&
\bigg(
\begin{array}{c}
\xi_{1}\\
\xi_{p}
\end{array}
\bigg)
+
\frac{2\tau}{p+2}
\bigg(
\begin{array}{c}
p-1 \\
-1
\end{array}
\bigg)
=:
\bigg(
\begin{array}{c}
\eta_1 \\
\eta_p
\end{array}
\bigg)
.
\end{eqnarray*}
Clearly, 
$
\eta_1
$
is the largest root of the polynomial
${\rm det}({\Zb}+ {\rm diag}(0, -v_\tau, \ldots, -v_\tau) + (2(p-1)\tau/(p+2)) {\bf I}_p - h {\bf I}_p)$, that is the largest eigenvalue of~${\Zb}+(2\tau/(p+2)){\rm diag}(p-1,-1,\ldots,-1)$, 
whereas 
$
\eta_p
$
is the smallest root of the polynomial~${\rm det}({\Zb}+ {\rm diag}(v_\tau,0, \ldots, 0) - (2\tau/(p+2)) {\bf I}_p - h {\bf I}_p)$, that is, the smallest eigenvalue of~${\Zb}+(2\tau/(p+2)){\rm diag}(p-1,-1,\ldots,-1)$. This proves the result for~$\tau\geq 0$.

We turn to Case~(b), for which~$\Sigb_n$ has eigenvalues 
\begin{equation}
	\label{eigenvic2}
\lambda_{n1}
=
\ldots
=
\lambda_{n,p-1}
=
\frac{1}{p} - \frac{2\tau}{\sqrt{n}p(p+2)}
\quad
\textrm{ and }
\quad
\lambda_{np}
=
\frac{1}{p} + \frac{2(p-1)\tau}{\sqrt{n}p(p+2)}
\cdot
\end{equation}
Here, we accordingly fix arbitrarily~$\thetab_1, \ldots, \thetab_{p-1}$ such that the $p \times p$ matrix~$\tilde{{\bf G}}_p:=(\thetab_1, \ldots, \thetab_{p-1},\thetab)$ is orthogonal. Still with~$\Lamb_n:={\rm diag}(\lambda_{n1},\lambda_{n2}, \ldots, \lambda_{np})$, we have $\Sigb_n \tilde{{\bf G}}_p=\tilde{{\bf G}}_p\Lamb_n$, so that $\tilde{{\bf G}}_p$ is an eigenvectors matrix for $\Sigb_n$. As above, $\xi_{n1}:=\sqrt{n}p (\hat{\lambda}_{n1}-\lambda_{n1})$ is the largest root of the polynomial~$P_{n1}(h):={\rm det}({\sqrt n}p({\bf S}_n-\lambda_{n1} {\bf I}_p)- h {\bf I}_p)$, whereas~$\xi_{np}:=\sqrt{n}p (\hat{\lambda}_{np}-\lambda_{np})$ is the smallest root of~$P_{np}(h):={\rm det}({\sqrt n}p({\bf S}_n-\lambda_{np} {\bf I}_p)- h {\bf I}_p)$. Letting
\begin{equation}
	\label{defZn2}
{\bf Z}_n
:= 
\sqrt{n}p(\tilde{{\bf G}}_p\pr {\bf S}_n \tilde{{\bf G}}_p - \Lamb_n)
=
 \tilde{{\bf G}}_p\pr \sqrt{n}p({\bf S}_n- \Sigb_n) \tilde{{\bf G}}_p
 ,
\end{equation}
rewrite these polynomials as
$
P_{nj}(h)
=
{\rm det}(\sqrt{n}p\tilde{{\bf G}}_p\pr {\bf S}_n \tilde{{\bf G}}_p - \sqrt{n}p\lambda_{nj} {\bf I}_p - h {\bf I}_p) 
=
{\rm det}({\Zb}_n+\sqrt{n}p(\Lamb_n-\lambda_{nj} {\bf I}_p)- h {\bf I}_p)
,
$
$j=1,p$, which here gives
$$
P_{n1}(h)
=
{\rm det}({\Zb}_n+{\rm diag}(0, \ldots,0, v_\tau)- h {\bf I}_p)
\nonumber  
$$
and
$$
P_{np}(h)
=
{\rm det}({\Zb}_n+{\rm diag}(-v_\tau,\ldots,-v_\tau,0)- h {\bf I}_p)
\nonumber  
,
$$
with the same~$v_\tau$ as above. Of course, we still have that~${\rm vec}\,\Zb_n
=
\sqrt{n}p 
(\tilde{{\bf G}}_p\otimes \tilde{{\bf G}}_p)' 
\, 
{\rm vec}\big(
\Sb_n
-
\Sigb_n
 \big)
$ converges weakly to~${\rm vec}\,\Zb\sim\mathcal{N}
(
{\bf 0},
{\bf V}_p
)
$. It readily follows that $(\xi_{n1},\xi_{np})'$ converges weakly to~$(\xi_{1},\xi_{p})'$, where~$\xi_{1}$ is the largest root of the polynomial~${\rm det}({\Zb}+ {\rm diag}(0, \ldots, 0,v_\tau) - h {\bf I}_p)$ and~$\xi_{p}$ is the smallest root of the polynomial~${\rm det}({\Zb}+{\rm diag}(-v_\tau, \ldots, -v_\tau,0)- h {\bf I}_p)$. This still implies that 
\begin{eqnarray*}
\bigg(
\begin{array}{c}
\sqrt{n} (p\hat{\lambda}_{n1}-1) \\
\sqrt{n} (p\hat{\lambda}_{np}-1) 
\end{array}
\bigg)
&\!\!=\!\!&
\bigg(
\begin{array}{c}
\xi_{n1}\\
\xi_{np}
\end{array}
\bigg)
+
\frac{2\tau}{p+2}
\bigg(
\begin{array}{c}
-1 \\
p-1
\end{array}
\bigg)
\\[2mm]
&\!\!\stackrel{\mathcal{D}}{\to}\!\!&
\bigg(
\begin{array}{c}
\xi_{1}\\
\xi_{p}
\end{array}
\bigg)
+
\frac{2\tau}{p+2}
\bigg(
\begin{array}{c}
-1 \\
p-1
\end{array}
\bigg)
=:
\bigg(
\begin{array}{c}
\eta_1 \\
\eta_p
\end{array}
\bigg)
.
\end{eqnarray*}
Clearly, 
$
\eta_1
$
is the largest root of the polynomial
${\rm det}({\Zb}+ {\rm diag}(0,\ldots,0,v_\tau) - (2\tau/(p+2)) {\bf I}_p - h {\bf I}_p)$, that is the largest eigenvalue of~${\Zb}+(2\tau/(p+2)){\rm diag}(-1,\ldots,-1,p-1)$, 
whereas 
$
\eta_p
$
is the smallest root of the polynomial~${\rm det}({\Zb}+ {\rm diag}(-v_\tau, \ldots,-v_\tau, 0) + (2(p-1)\tau/(p+2)) {\bf I}_p - h {\bf I}_p)$, that is, the smallest eigenvalue of~${\Zb}+(2\tau/(p+2)){\rm diag}(-1,\ldots,-1,p-1)$. This establishes the result for~$\tau<0$. 
\cqfd
\vspace{3mm}


{\sc Proof of Corollary~\ref{Corollambda1asympt}.}
According to Theorem~\ref{thlambda1asympt}, $L_p^{\rm max}$ is equal in distribution to the first marginal of~${\pmb \ell}^{(p)}=(\ell^{(p)}_{1},\ldots,\ell^{(p)}_{p})'$, where $\ell^{(p)}_{1}\geq \ldots \geq \ell^{(p)}_{p}$ are the eigenvalues of~$\Zb=(Z_{ij})$, with~${\rm vec}\,\Zb\sim\mathcal{N}({\bf 0},{\bf V}_p)$. Note that $Z_{ij}=Z_{ji}$ almost surely for any~$1\leq i<j\leq p$, so~${\rm vec}\,\Zb$ may not have a density with respect to the Lebesgue measure on~$\R^{p^2}$, but~${\rm vech}\,\Zb$ might in principle have a density with respect to the Lebesgue measure on~$\R^{p(p+1)/2}$. If~${\rm vech}\,\Zb$ indeed has such a density, then Theorem~13.3.1 of \cite{And2003} can be used to obtain the density of~${\pmb \ell}^{(p)}$ from that of~${\rm vech}\,\Zb$. However, since~${\rm tr}[\Zb]=({\rm vec}\,{\bf I}_p)'({\rm vec}\,\Zb)=0$ 
 \vspace{-.7mm}
almost surely, ${\rm vech}\,\Zb$ does not admit a density and the eigenvalues~$\ell^{(p)}_{1},\ldots,\ell^{(p)}_{p}$ add up to zero almost surely (and thus do not admit a joint density over~$\R^p$ either). 

We solve this issue by considering a sequence of $p\times p$ random matrices~$\Zb_{\delta_k}$ ($k=1,2,\ldots$), with~$\delta_k(> 0)$ converging to zero as~$k$ goes to infinity, and such that~${\rm vec}\,\Zb_\delta\sim\mathcal{N}
(
{\bf 0},
{\bf V}_{p,\delta}
)$, with
$$
{\bf V}_{p,\delta}
:=
\frac{p}{p+2} (\mathbf{I}_{p^2}+{\bf K}_p) - \frac{2-\delta}{p+2} \mathbf{J}_{p}
.
$$
Of course, $\Zb_{\delta_k}$ converges weakly to $\Zb$, so that the continuous mapping theorem ensures
 \vspace{-.4mm}
 that~${\pmb \ell}^{(p)}_{\delta_k}=(\ell^{(p)}_{1\delta_k},\ldots,\ell^{(p)}_{p\delta_k})'$ (here, $\ell^{(p)}_{1\delta}\geq \ldots \geq \ell^{(p)}_{p\delta}$ are the eigenvalues of~$\Zb_\delta$) also converges weakly to~${\pmb \ell}^{(p)}=(\ell^{(p)}_{1},\ldots,\ell^{(p)}_{p})'$.  
Let then~$\Db_p$ be the $p$-dimensional duplication matrix, that is such that~$\Db_p({\rm vech}\,\Ab)={\rm vec}\,\Ab$ for any~$p\times p$ symmetric matrix~$\Ab$. Write~$\Wb_\delta:={\rm vech}\,\Zb_\delta=\Db_p^-({\rm vec}\,\Zb_\delta)$, where~$\Db_p^-=(\Db_p'\Db_p)^{-1}\Db_p'$ is the Moore-Penrose inverse of~$\Db_p$. By definition of~$\Zb_\delta$, the random vector~$\Wb_\delta$ has density
$$ 
{\bf w}\mapsto  
h_\delta({\bf w}) 
=  
\frac{a_p}{\sqrt{\delta}}  
\exp\Big(-\frac{1}{2}\, {\bf w}' 
\big(
\Db_p^-
\big\{
{\textstyle{\frac{p}{p+2}}} 
({\bf I}_{p^2}+ {\bf K}_p)
-
{\textstyle{\frac{2-\delta}{p+2}}} 
 \Jb_p
\big\}
(\Db_p^-)'
\big)^{-1} 
{\bf w} \Big)
,
$$
where
$$
a_p
:=
\frac{(p+2)^{p(p+1)/4}}{2^{(p^2+3p-2)/4} (\pi p)^{p(p+1)/4}}
$$ 
is a normalizing constant. By using the identities~${\bf K}_p{\bf D}_p={\bf D}_p$ and $\Db_p^-({\rm vec}\,\mathbf{I}_p)=\Db_p'({\rm vec}\,\mathbf{I}_p)={\rm vech}\,\mathbf{I}_p$, it is easy to check that 
\begin{eqnarray*}
\lefteqn{
\big
(\Db_p^-
\big\{
{\textstyle{\frac{p}{p+2}}} 
({\bf I}_{p^2}+ {\bf K}_p)
-
{\textstyle{\frac{2-\delta}{p+2}}} 
 \Jb_p
\big\}
(\Db_p^-)'
\big)^{-1}
}
\\[2mm]
& & 
\hspace{3mm} 
=
\frac{p+2}{2p}\, 
\Db_p'
\bigg\{
\frac{1}{2} ({\bf I}_{p^2}+ {\bf K}_p)
+
\frac{\frac{2-\delta}{p+2}}{\frac{2p}{p+2}-\frac{p(2-\delta)}{p+2}} 
\Jb_p
\bigg\}
\Db_p
\\[2mm]
& & 
\hspace{3mm} 
=
\frac{p+2}{2p}\, 
\Db_p'
\bigg\{
\frac{1}{2} ({\bf I}_{p^2}+ {\bf K}_p)
+
\frac{2-\delta}{\delta p} 
\Jb_p
\bigg\}
\Db_p
.
\end{eqnarray*}
 Using the identities~$({\rm vec}\,\Ab)'({\rm vec}\,\Bb)={\rm tr}[\Ab'\Bb]$ and~${\bf K}_p({\rm vec}\,\Ab)={\rm vec}(\Ab')$, the resulting density for~$\Zb_\delta$ is therefore
\begin{eqnarray*}
\lefteqn{
\zb \mapsto f(\zb)=h_\delta({\rm vech}\,\zb)
}
\\[2mm]
& & 
\hspace{3mm} 
=
\frac{a_p}{\sqrt{\delta}}
\exp\Big(
-
\frac{p+2}{4p}\, 
({\rm vech}\,\zb)' 
\Db_p'
\bigg\{
\frac{1}{2} ({\bf I}_{p^2}+ {\bf K}_p)
+
\frac{2-\delta}{\delta p} 
\Jb_p
\bigg\}
\Db_p
({\rm vech}\,\zb)
 \Big)
\\[2mm]
& & 
\hspace{3mm} 
=
\frac{a_p}{\sqrt{\delta}}
\exp\Big(
-
\frac{p+2}{4p}\, 
\, 
({\rm vec}\,\zb)'
\Big\{
{\bf I}_{p^2}
+
\frac{2-\delta}{\delta p} 
\Jb_p
\Big\}
{\rm vec}\,\zb \Big)
\\[2mm]
& & 
\hspace{3mm} 
=
\frac{a_p}{\sqrt{\delta}}
\exp\Big(
-
\frac{p+2}{4p}\, 
\, 
\Big\{
({\rm tr}\,\zb^2)
+
\frac{2-\delta}{\delta p} 
({\rm tr}\,\zb)^2
\Big\}
\Big)
.
\end{eqnarray*}
Theorem~13.3.1 from \cite{And2003} then implies that~${\pmb \ell}^{(p)}_{\delta}=(\ell^{(p)}_{1\delta},\ldots,\ell^{(p)}_{p\delta})'$ has density
\begin{eqnarray*}
	\lefteqn{
\hspace{-10mm} 
(\ell_1, \ldots, \ell_p)' 
\mapsto 
\frac{b_p}{\sqrt{\delta}}
\exp
\bigg(
-
\frac{p+2}{4p}
\, 
\bigg\{
\bigg( \sum_{j=1}^p \ell_j^2 \bigg)
+
\frac{2-\delta}{\delta p} 
\bigg(\sum_{j=1}^p \ell_j\bigg)^2
\bigg\}
\bigg)
}
\\[2mm]
& & 
\hspace{40mm} 
\times
\bigg( \prod_{1\leq k<j \leq p} (\ell_k- \ell_j) \bigg)
\,
\mathbb{I}[\ell_1\geq\ldots\geq \ell_p]
,
\nonumber
\end{eqnarray*}
with
$$
b_p
:=
\frac{(p+2)^{p(p+1)/4}}{2^{(p^2+3p-2)/4} p^{p(p+1)/4}\prod_{j=1}^{p} \Gamma(\frac{j}{2})}
\cdot
$$ 

We now turn to the particular cases (i)~$p=2$ and~(ii)~$p=3$ considered in the statement of the corollary. (i) Our general derivation above states that ${\pmb \ell}^{(2)}_{\delta}=(\ell^{(2)}_{1\delta},\ell^{(2)}_{2\delta})'$ has density~$(\ell_1,\ell_2)' 
\mapsto \mathcal{I}(\ell_1,\ell_2)
\mathbb{I}[\ell_1\geq\ell_2]$, with
$$
\mathcal{I}(\ell_1,\ell_2)
:=
\frac{1}{\sqrt{2\pi\delta}}
(\ell_1-\ell_2) 
e^{
-
\frac{1}{2}
\{
( \ell_1^2+\ell_2^2 )
+
\frac{2-\delta}{2\delta} 
(\ell_1+\ell_2)^2
\}
}
.
$$
Direct computations allow checking that~$\mathcal{I}(\ell_1,\ell_2)$ is the derivative of the function
$$
\ell_2
\mapsto
\frac{\sqrt{2\delta}}{\sqrt{\pi}(2+\delta)}
e^{
-
\frac{1}{2}
\{
( \ell_1^2+\ell_2^2 )
+
\frac{2-\delta}{2\delta} 
(\ell_1+\ell_2)^2
\}
}
 +
\frac{2^{5/2} \ell_1}{(2+\delta)^{3/2}}
e^{
-\frac{2\ell^2_1}{2+\delta}
}
\Phi
\bigg(
\frac{(2-\delta)\ell_1+(2+\delta)\ell_2}{\sqrt{2\delta(\delta+2)}}
\bigg)
.
$$
It follows that~$\ell^{(2)}_{1\delta}$ has density 
\begin{eqnarray*}
\lefteqn{	
\hspace{0mm} 
\ell_1 
\mapsto
\int_{-\infty}^{\ell_1} 
\mathcal{I}(\ell_1,\ell_2)
\,d\ell_2
}
\\[2mm]
 & &
 \hspace{-4mm} 
=
\bigg[
\frac{\sqrt{2\delta}}{\sqrt{\pi}(2+\delta)}
e^{
-
\frac{1}{2}
\{
( \ell_1^2+\ell_2^2 )
+
\frac{2-\delta}{2\delta} 
(\ell_1+\ell_2)^2
\}
}
 +
\frac{2^{5/2} \ell_1}{(2+\delta)^{3/2}}
e^{
-\frac{2\ell^2_1}{2+\delta}
}
\Phi
\bigg(
\frac{(2-\delta)\ell_1+(2+\delta)\ell_2}{\sqrt{2\delta(\delta+2)}}
\bigg)
\bigg]_{-\infty}^{\ell_1}
\\[2mm]
 & &
 \hspace{-4mm} 
=
\frac{\sqrt{2\delta}}{\sqrt{\pi}(2+\delta)}
e^{
-\frac{2\ell^2_1}{\delta}
}
+
\frac{2^{5/2} \ell_1}{(2+\delta)^{3/2}}
e^{
-\frac{2\ell^2_1}{2+\delta}
}
\Phi
\bigg(
\frac{4\ell_1}{\sqrt{2\delta(\delta+2)}}
\bigg)
.
\end{eqnarray*}
Taking the limit as~$\delta\to 0$, we obtain that~$\ell^{(2)}_{1}(\stackrel{\mathcal{D}}{=}L_2^{\rm max})$ has density
$
\ell_1
\mapsto
2\ell_1
\exp(-\ell^2_1)
\mathbb{I}[\ell_1>0]
$. From Theorem~\ref{thlambda1asympt}, this is the density of the asymptotic distribution of~$T_+\n=\sqrt{n} (p\hat{\lambda}_{n1}-1)$. Part~(i) of the result then follows from the fact that~$T_+\n=T_-\n=T_\pm\n$ almost surely for~$p=2$ (see the discussion below the corollary). 

(ii) For~$p=3$, the density of~${\pmb \ell}^{(3)}_{\delta}=(\ell^{(3)}_{1\delta},\ell^{(3)}_{2\delta},\ell^{(3)}_{3\delta})'$ is~$(\ell_1,\ell_2,\ell_3)' \mapsto \mathcal{J}(\ell_1,\ell_2,\ell_3)
\mathbb{I}[\ell_1\geq\ell_2\geq\ell_3]$, with
$$
\mathcal{J}(\ell_1,\ell_2,\ell_3)
:=
\frac{125}{216\pi\sqrt{\delta}}
(\ell_1-\ell_2)
(\ell_1-\ell_3)
(\ell_2-\ell_3) 
e^{
-
\frac{5}{12}
\, 
\{
( \ell_1^2+\ell_2^2+\ell_3^2 )
+
\frac{2-\delta}{3\delta} 
(\ell_1+\ell_2+\ell_3)^2
\}
}
.
$$
Lengthy yet straightforward computations show that~$\mathcal{J}(\ell_1,\ell_2,\ell_3)$ is the derivative of the function~$\ell_3
\mapsto \mathcal{K}(\ell_1,\ell_2,\ell_3)$, with
\begin{eqnarray*}
\lefteqn{	
\hspace{0mm} 
\mathcal{K}(\ell_1,\ell_2,\ell_3)
:=
}
\\[2mm]
 & &
\hspace{-4mm} 
-
\frac{25\sqrt{\delta}(\ell_1-\ell_3)}{48\pi(1+\delta)^2}
\{2(2\ell_1+2\ell_3-\ell_2)+\delta(\ell_1+\ell_3-2\ell_2)\}
e^{
-
\frac{5}{12}
\, 
\{
( \ell_1^2+\ell_2^2+\ell_3^2 )
+
\frac{2-\delta}{3\delta} 
(\ell_1+\ell_2+\ell_3)^2
\}
}
\\[2mm]
 & &
\hspace{-4mm} 
-
\frac{5\sqrt{10}(\ell_1-\ell_3)}{288\sqrt{\pi}(1+\delta)^{5/2}}
\{20(2\ell^1_1+5\ell_1\ell_3+2\ell_3^2)-2\delta(5(\ell_1-\ell_3)^2-18)-\delta^2(5(\ell_1-\ell_3)^2-36)\}
\\[2mm]
 & &
 \hspace{10mm} 
\times
e^{
-
\frac{20(\ell_1^2+\ell_1\ell_3+\ell_3^2)+5\delta (\ell_1-\ell_3)^2)}{24(1+\delta)}
}
\Phi
\bigg(
\frac{\sqrt{5}(2(\ell_1+\ell_2+\ell_3)-\delta(\ell_1+\ell_3-2\ell_2))}{6\sqrt{\delta(1+\delta)}}
\bigg)
.
\end{eqnarray*}
Therefore, $(\ell^{(3)}_{1\delta},\ell^{(3)}_{3\delta})'$ has density 
$$
(\ell_1,\ell_3) 
\mapsto
\bigg(
\int_{\ell_3}^{\ell_1} 
\mathcal{J}(\ell_1,\ell_2,\ell_3)
\,d\ell_2
\bigg)
\mathbb{I}[\ell_1\geq\ell_3]
=
(\mathcal{K}(\ell_1,\ell_1,\ell_3)-\mathcal{K}(\ell_1,\ell_3,\ell_3))
\mathbb{I}[\ell_1\geq\ell_3]
,
$$
with
\begin{eqnarray*}
\lefteqn{	
\hspace{-2mm} 
\mathcal{K}(\ell_1,\ell_1,\ell_3)-\mathcal{K}(\ell_1,\ell_3,\ell_3)
=
}
\\[2mm]
 & &
\hspace{-3mm} 
-
\frac{25\sqrt{\delta}(\ell_1-\ell_3)}{48\pi(1+\delta)^2}
\{2(\ell_1+2\ell_3)+\delta(\ell_3-\ell_1)\}
e^{
-
\frac{5}{12}
\, 
\{
( 2\ell_1^2+\ell_3^2 )
+
\frac{2-\delta}{3\delta} 
(2\ell_1+\ell_3)^2
}
\\[2mm]
 & &
\hspace{-3mm} 
+
\frac{25\sqrt{\delta}(\ell_1-\ell_3)}{48\pi(1+\delta)^2}
\{2(2\ell_1+\ell_3)+\delta(\ell_1-\ell_3)\}
e^{
-
\frac{5}{12}
\, 
\{
( \ell_1^2+2\ell_3^2 )
+
\frac{2-\delta}{3\delta} 
(\ell_1+2\ell_3)^2
}
\\[2mm]
 & &
\hspace{-3mm} 
-
\frac{5\sqrt{10}(\ell_1-\ell_3)}{288\sqrt{\pi}(1+\delta)^{5/2}}
e^{
-
\frac{20(\ell_1^2+\ell_1\ell_3+\ell_3^2)+5\delta (\ell_1-\ell_3)^2)}{24(1+\delta)}
}
\\[2mm]
 & &
\hspace{4mm} 
\times
\{20(2\ell^1_1+5\ell_1\ell_3+2\ell_3^2)-2\delta(5(\ell_1-\ell_3)^2-18)-\delta^2(5(\ell_1-\ell_3)^2-36)\}
\\[2mm]
 & &
\hspace{4mm} 
\times
\bigg[
\Phi
\bigg(
\frac{\sqrt{5}(2(2\ell_1+\ell_3)-\delta(\ell_3-\ell_1))}{6\sqrt{\delta(1+\delta)}}
\bigg)
-
\Phi
\bigg(
\frac{\sqrt{5}(2(\ell_1+2\ell_3)-\delta(\ell_1-\ell_3))}{6\sqrt{\delta(1+\delta)}}
\bigg)
\bigg]
.
\end{eqnarray*}
Taking the limit as~$\delta\to 0$ shows that the density of~$(\ell^{(3)}_{1},\ell^{(3)}_{3})'$ is 
\begin{equation}
\label{jointlambdap3}	
(\ell_1,\ell_3) 
\mapsto
-
\frac{100\sqrt{10}}{288\sqrt{\pi}}
(\ell_1-\ell_3)(2\ell^2_1+5\ell_1\ell_3+2\ell_3^2)
e^{
-
\frac{5(\ell_1^2+\ell_1\ell_3+\ell_3^2)}{6}
}
\mathbb{I}[-2\ell_1< \ell_3\leq -\ell_1/2]
.
\end{equation}
From Theorem~\ref{thlambda1asympt}, the density of the asymptotic null distribution of~$T_+\n=\sqrt{n} (p\hat{\lambda}_{n1}-1)$ coincides with the density of~$\ell_1^{(3)}(\stackrel{\mathcal{D}}{=}L_3^{\rm max})$. After marginalization in~(\ref{jointlambdap3}), this last density is seen to be
\begin{eqnarray*}
\lefteqn{
\ell_1
\mapsto
\bigg\{
\sqrt{\frac{5}{2\pi}} e^{-\frac{5\ell_1^2}{2}}
+
\sqrt{\frac{5}{8\pi}} e^{-\frac{5\ell_1^2}{8}}
+
\frac{3}{4}\sqrt{\frac{5}{8\pi}} (5\ell_1^2-4) e^{-\frac{5\ell_1^2}{8}}
\bigg\}
\mathbb{I}[\ell_1\geq 0]
}
\\[2mm]
 & & 
 \hspace{3mm} 
=
\bigg\{
\frac{d}{d\ell_1} \Phi(\sqrt{5}\ell_1)
+
\frac{d}{d\ell_1} \Phi({\textstyle{\frac{\sqrt{5}\ell_1}{2}}})
+
3\frac{d}{d\ell_1} \Phi''({\textstyle{\frac{\sqrt{5}\ell_1}{2}}})
\bigg\}
\mathbb{I}[\ell_1\geq 0]
,
 \end{eqnarray*}
which proves the result for~$T_+\n$, hence also for~$T_-\n$ (recall from the discussion below the corollary that~$T_+\n$ and~$T_-\n$ share the same weak limit in any dimension~$p$). Finally, the result for~$T_\pm\n$ follows from~(\ref{jointlambdap3}) by using 
the fact that~$T_\pm\n$ converges weakly to~$\max(\ell_1^{(3)},-\ell_3^{(3)})$.
\cqfd
\vspace{3mm}

\section*{Acknowledgement}

Davy Paindaveine's research is supported by a research fellowship from the Francqui Foundation and by the Program of Concerted Research Actions (ARC) of the Universit\'{e} libre de Bruxelles. Thomas Verdebout's research is supported by the ARC Program of the Universit\'{e} libre de Bruxelles and by the Cr\'{e}dit de Recherche J.0134.18 of the FNRS (Fonds National pour la Recherche Scientifique), Communaut\'{e} Fran\c{c}aise de Belgique.


\bibliographystyle{imsart-nameyear}
\bibliography{Paper}           
\vspace{3mm}

\end{document}